# Selmer groups and the Eisenstein-Klingen Ideal

E.Urban*

March 24, 1998

## 0   Introduction

The central point in the Bloch-Kato conjectures is to establish formulas for the order of the Selmer groups attached to Galois representations in terms of the special values of their L-functions. In order to give upper bound, the main way is to construct Euler systems following Kolyvagin. Besides, lower bounds have been obtained by using congruences between automorphic forms. So, a lot of Iwasawa conjectures for number fields have been attacked or proven by using congruences between modular forms for $GL(2)$ over the corresponding number field (cf. [22], [36], [21], [18] and [30]). It seems that the new developments of the arithmetic theory of automorphic forms for bigger groups than $GL_2$ enable us to apply again this idea for non abelian Iwasawa theory.

In this circle of ideas, this article is devoted to the first part of a strategy of proving a divisibility towards a Main Conjecture "à la Iwasawa" as it is announced in [19] for the Selmer group of a two variable adjoint modular Galois representation by using congruences between cuspidal Siegel modular forms of genus 2 and Klingen type Eisenstein series.

In order to be more precise in the formulation of that conjecture, let us introduce some notations. Let $\mathbf{I}$ be a finite and flat extension of $\mathbf{Z}_p[[X]]$. Let $\mathcal{F} = \sum_{n\geq 1} a(n;\mathcal{F})q^n \in \mathbf{I}[[q]]$ be an $\mathbf{I}$-adic cuspidal modular forms of tame level $N$ and nebentypus $\epsilon_\mathcal{F}$ (i.e. $\mathcal{F}$ modulo $P$ is an eigen $p$-ordinary cuspidal modular form of weight $k$, level $Np$ and nebentypus $\epsilon_\mathcal{F}\omega^{-k}$ if $P$ is a height one prime ideal of $\mathbf{I}$ above $(1+T-u^k)$ for $u$ a topological generator of $1+p\mathbf{Z}_p$ fixed once for all and $\omega$ the Teichmüller character). It is well known that (when $\mathcal{F}$ is not residually "Eisenstein") one can associate to $\mathcal{F}$ a Galois representation $\rho_{\mathcal{F}} : G_{\mathbf{Q}} \to GL_2(\mathbf{I})$. Let $\eta$ be a Dirichlet character of level $Np$; then one can consider $\tilde{\eta} : G_{\mathbf{Q}} \to \mathbf{Z}_p[[S]]^{\times}$ the universal deformation of the Galois character associated to $\eta$. On one hand, one can define a Selmer group associated to $ad^0(\rho_{\mathcal{F}}) \otimes \tilde{\eta}^{-1}$ (see paragraph 4.). Generalizing Iwasawa conjecture, Greenberg conjectures in [8] that the Selmer group $Sel(ad^0(\rho_{\mathcal{F}}) \otimes \tilde{\eta}^{-1})$ is co-torsion over $\mathbf{I}[[S]]$; this has been proven by Hida in the case $\eta = 1$ cf. [15] and his proof should extend (at least in some cases cf. [16]) for even characters $\eta$ thanks to Fujiwara's generalization of the results of Taylor and Wiles (cf. [7]). We assume this conjecture for a moment and denote by $F_{ad^0(\rho_{\mathcal{F}}) \otimes \tilde{\eta}^{-1}}$ the characteristic ideal of this Selmer group.

---

*CNRS-LAGA Université de Paris-Nord, France and UCLA Department of Mathematics, Los Angeles, CA 90095-1555 U.S.A. This work was partially supported by an NSF grant DMS 97701017.

On the other hand, almost ten years ago, in [14] Hida constructed a two-variable $p$-adic L-function associated to $\mathcal{F}$ and $\eta$ (when $\eta$ is even), say $L \in Frac(\mathbf{I}[[S]])$, verifying the following interpolation property: For $P$ above $(1 + T - u^k)$ and $Q = (1 + S - u^l)$ satisfying $-k \leq l \leq 0$, we have

$$L(P,Q) = \star E(k,l) \frac{L(1-l, Ad(\rho_{f_k}) \otimes \eta^{-1}\omega^l)}{(2i\pi)^{-2l}\Omega(f_k)}$$

for a constant $\star$, $E(k,l)$ a simple Euler factors at $p$ and $\Omega(f_k) = (2i)^{k+1}\pi^2 < f_k, f_k >$ with $f_k$ the newform associated to $\mathcal{F} \bmod P$. If $H$ annihilates the congruence module associated to $\mathcal{F}$, it is proven in [14] that $H \times L \in \mathbf{I}[[S]]$. Moreover if $\eta_p$ (the $p$-component of $\eta$) is trivial $ad^0(\rho_{\mathcal{F}})$ and $\epsilon_{\mathcal{F}} = \omega^2$, there is a trivial zero in $(S)$ since $E(k,l)$ vanishes on the line $l = 0$. We denote by $L_{ad^0(\rho_{\mathcal{F}}) \otimes \tilde{\eta}^{-1}} \in \mathbf{I}[[S]]$ the product of $L$ and the characteristic ideal of the congruence module of $\mathcal{F}$ ( that we divided by $S$ in the case $\eta_p$ trivial and $\epsilon_{\mathcal{F}} = \omega^2$). Any generator interpolates up to a p-adic unit the critical values divided $(2i\pi)^{-2l}\Omega_{Hida}(f_k)$ where $\Omega_{Hida}(f_k)$ is the Hida's normalization (cf. [11]) of the Deligne's period $\Omega(f_k)$ associated to $f_k$.

Let $p^* = p(-1)^{\frac{p-1}{2}}$. the Main conjecture formulated by Greenberg in [8] asserts:

**Conjecture 0.1** *Assume $\bar{\rho}_{\mathcal{F}}$ restricted to $Gal(\bar{\mathbf{Q}}/\mathbf{Q}(\sqrt{p^*}))$ is absolutely irreducible, then we have the equality*

$$L_{ad^0(\rho_{\mathcal{F}}) \otimes \tilde{\eta}^{-1}} = F_{ad^0(\rho_{\mathcal{F}}) \otimes \tilde{\eta}^{-1}}$$

In this article, we are setting up a strategy to prove one divisibility towards that conjecture. The idea is to introduce a third characteristic ideal we call $Eis(\mathcal{F}, \eta)$ containing informations on the congruences between cuspidal Siegel modular forms and the Eisenstein series constructed "inducing" $\mathcal{F} \otimes \tilde{\eta}$ from the Klingen parabolic subgroup to $GSp_4$ and to prove the two divisibilities:

(i) $Eis(\mathcal{F}, \eta) \mid F_{ad^0(\rho_{\mathcal{F}}) \otimes \tilde{\eta}^{-1}}$

(ii) $L_{ad^0(\rho_{\mathcal{F}}) \otimes \tilde{\eta}^{-1}} \mid Eis(\mathcal{F}, \eta)$

In [19], it is explained how one can prove the above conjecture if $\eta_p$ is trivial and $\epsilon_{\mathcal{F}} = \omega^2$) from the divisibility $L_{ad^0(\rho_{\mathcal{F}}) \otimes \tilde{\eta}^{-1}} \mid F_{ad^0(\rho_{\mathcal{F}}) \otimes \tilde{\eta}^{-1}}$ by using the results of Hida [15] and Greenberg-Tilouine [9]. Let us now state the main result of this paper.

**Theorem 0.1** *Under some hypothesis (see theorem 3.5 for a precise statement), the following divisibility holds:*

$$Eis(\mathcal{F}, \eta) \mid F_{ad^0(\rho_{\mathcal{F}}) \otimes \tilde{\eta}^{-1}}$$

Actually, we give a more general result including the case of one variable Selmer group attached to the cyclotomic deformations of the adjoint Galois representation associated to a modular form. The second divisibility is the subject of a subsequent paper (cf. [32]) and is motivated by the fact that constant terms of the Klingen type Eisenstein series are expressed in terms of the critical values interpolated by $L_{ad^0(\rho_{\mathcal{F}}) \otimes \tilde{\eta}^{-1}}$.

As written above, we make use of congruences between cuspidal Siegel modular forms and Klingen-type Eisenstein series in order to construct sufficiently indecomposable Galois modules with suitable semi-simplification. The Eisenstein ideal is constructed by the usual way as an ideal of the universal ordinary Hecke algebra for $GSp_4$ whose properties have been developed by J. Tilouine and the author in [28] (cf. section 2).

The construction of the non trivial extension uses the existence of the 4-dimensional Galois representations associated to the cuspidal representations of $GSp_4(\mathbf{A_A})$ recently proven by R. Weissauer. The point is that the Galois representations associated to the Eisenstein series are reducible while those associated to cuspidal representation (if they are neither CAP nor Endoscopic) have to be absolutely irreducible. This last point is still unproven in general. However, we have been able to prove it (see theorem 3.3) in our situation thanks to the Taylor-Wiles'theorem (actually Diamond's improvement cf. [4]) asserting that an ordinary deformation of a modular Galois representation is modular.

Once we know the irreducibility of these Galois representations, the second step is to use the congruence property to construct non split exact sequences of type

$$0 \to \rho_{\mathcal{F}} \otimes \mathbf{I}[[S]]_P / P^r \to M \to \rho_{\mathcal{F}} \otimes \tilde{\eta} \otimes \mathbf{I}[[S]]_P / P^r \to 0$$

Then the non-splitness enable us to construct cocycles in our Selmer group. This is what Mazur and Wiles did in their proof of the main Iwasawa conjecture for $\mathbf{Q}$ (and the Wiles'generalization to totally real fields). However, their method used the abelian context and could not be generalized here. Therefore we needed to produce these extensions in a totally different way (see section 1).

The expected properties of the constructed cocycles result from the properties of the Galois representations we use. On the one hand, we assume temporarily that the $p$-adic Galois representations associated to $p$-ordinary cuspidal representations are ordinary (see paragraph 3 for a discussion of what it is know for the local properties at $p$ of these Galois representations). On the other hand, we ask that they take values in $GSp_4(\bar{\mathbf{Q}}_p)$ in order to get a cocycle with trace zero. By our irreducibility result, it is easy to see these representations respect a skew-symmetric or a symmetric bilinear form, but we do not know yet how to exclude this last case in general. A discussion of that question is made at the end of paragraph 3.2.

It seems that this approach can be developed in many situations. As an other example, one can expect seriously that considering the group $GL(3)$ one can get similarly results for the Selmer group of the standard Galois representation associated to ordinary Hilbert modular forms considering cuspidal representations of $GL(3)_{/F}$ congruent to Eisenstein series associated to a maximal parabolic. We note also that the results of the first part of this paper suggests that a general theory of Eisenstein Ideal and of deformations of reducible Galois representation can be developed in order to give systematically lower bound for more general Selmer groups. These considerations are developed in [33] and [34].

**Acknowledgments:** A part of this paper was worked out during my visit to the Mehta Research Institute (MRI) of Allahabad, so I take the opportunity to thank D. Prasad for his kind invitation. It is for me a pleasure to thank warmly H. Hida and J. Tilouine for helpful conversations during the preparation of this work.

# Contents

# 1 Residually reducible representations on local rings

## 1.1 Residually reducible representations on lattices.

We now give now a variant of the result of [31] we need in the last section of this article. For that purpose, we introduce the following notations:

Let $\mathcal{B}$ be a henselian, generically etale and reduced local commutative algebra which is finite over a discrete valuation ring $\mathcal{O}$. We denote by $\mathcal{M}_\mathcal{B}$ (resp. $\kappa_\mathcal{B}$, $\tilde{\mathcal{B}}$ and $F_\mathcal{B}$ ) the maximal ideal of $\mathcal{B}$ (resp. the residue field, the normalization and the total ring of fractions of $\mathcal{B}$). Since $\mathcal{B}$ is reduced, we can embed it in the product of its irreducible components:

$$\mathcal{B} \subset \tilde{\mathcal{B}} = \prod_x \mathcal{O}_x \subset \prod_x F_x = F_\mathcal{B}$$

with $F_x$ the field of fractions of the irreducible component $\mathcal{O}_x$.

**Theorem 1.1** *Let $\mathcal{R}$ be a $\mathcal{B}$-algebra and $\rho$ be an absolutely irreducible representation of $\mathcal{R}$ on $F_\mathcal{B}^n$ such that there exist two representations $\rho_i$ for $i = 1, 2$ in $M_{n_i}(\mathcal{B})$ and $I \subset \mathcal{B}$ a proper ideal of $\mathcal{B}$ such that:*

(i) *The coefficients of the characteristic polynomial of $\rho$ belongs to $\mathcal{B}$.*

(ii) *The characteristic polynomials of $\rho$ and $\rho_1 \oplus \rho_2$ are congruent modulo $I$.*

(iii) *$\bar{\rho}_1$ and $\bar{\rho}_2$ are absolutely irreducible*

(iv) *$\bar{\rho}_1 \neq \bar{\rho}_2$*

*Then there exists a $\mathcal{R}$-stable $\mathcal{B}$-lattice $\mathcal{L}$ in $F_\mathcal{B}^n$ and a finitely generate $\mathcal{B}$-submodule $\mathcal{T}$ of $F_\mathcal{B}$ such that we have the following exact sequence:*

$$0 \longrightarrow \rho_1 \otimes \mathcal{T}/I.\mathcal{T} \longrightarrow \mathcal{L} \otimes \mathcal{B}/I \xrightarrow{s} \rho_2 \otimes \mathcal{B}/I \longrightarrow 0$$

*where $s$ is only a section of $\mathcal{B}/I$-module. Moreover $\mathcal{L}$ has no quotient isomorphic to $\bar{\rho}_1$.*

**Proof:** For each $x$, let us denote by $\mathcal{M}_x$ the maximal ideal of $\mathcal{O}_x$. Let $\rho_x$ be the representation on $F_x^n$ deduced from $\rho$. Of course, it is absolutely irreducible by our hypothesis. Let now $\mathcal{L}_x$ be a stable $\mathcal{O}_x$-lattice of $F_x^n$ (which exists since by (i) the trace of $\rho_x$ takes values in $\mathcal{O}_x$; this can be seen easily by arranging the proof of lemma 6 of [24]) and consider a minimal stable sublattice of $\mathcal{L}_x$ such that $\mathcal{L}_x \otimes \mathcal{O}_x/\mathcal{M}_x$ contains a stable subspace isomorphic to $\bar{\rho}_1$; the latter exists because characteristic polynomials of $\rho_x$ and $\rho_1 \oplus \rho_2$ are congruent modulo $\mathcal{M}_x$. Taking a suitable $\mathcal{O}_x$-basis of $\mathcal{L}_x$, we can get a representation:

$$\rho_x : \mathcal{R} \to M_n(\mathcal{O}_x) \text{ such that } \bar{\rho}_x(r) = \begin{pmatrix} \bar{\rho}_1(r) & \star \\ 0 & \bar{\rho}_2(r) \end{pmatrix}.$$

We then set $\rho_{\tilde{\mathcal{B}}} = (\rho_x)_x$ which takes values in $M_n(\tilde{\mathcal{B}})$ and

$$\bar{\rho}_{\tilde{\mathcal{B}}}(r) = \begin{pmatrix} \bar{\rho}_1(r) & \star \\ 0 & \bar{\rho}_2(r) \end{pmatrix} \in M_n(\tilde{\mathcal{B}}/Rad(\tilde{\mathcal{B}})).$$

If $r_0 \in \mathcal{R}$ is such that the characteristic polynomial of $\bar{\rho}(r_0)$ has exactly $n$ different roots in $\kappa_{\mathcal{B}}$, by the Hensel lemma the characteristic polynomial of $\rho_{\tilde{\mathcal{B}}}(r_0)$ has $n$ roots $\alpha_1, \ldots, \alpha_n$ in $\mathcal{B}$ which are different modulo the maximal ideal and thus we can assume (as in the proof of Theorem of [31]) that $\rho_{\tilde{\mathcal{B}}}(r_0) = Diag(\alpha_1, \ldots, \alpha_n)$ and the sub-$\mathcal{B}$-module of $M_n(\tilde{\mathcal{B}})$ generated by the powers of $\rho_{\tilde{\mathcal{B}}}(r_0)$ is exactly the set of diagonal matrices with entries in $\mathcal{B}$. For each $i$, $1 \leq i \leq n$, we can thus choose $s_i \in \mathcal{R}$ such that

$$\rho_{\tilde{\mathcal{B}}}(s_i) = E(\alpha_i) = Diag(0, \ldots, 0, 1, 0, \ldots, 0)$$

where 1 is at the $i$-th position. Of course, $E(\alpha_i)$ is the projector of the eigenspace of $r_0$ associated to the eigenvalue $\alpha_i$. There exist also elements $r_1, r_2 \in \mathcal{R}$ such that

$$\rho_{\tilde{\mathcal{B}}}(r_1) = E_1 = \begin{pmatrix} \mathbf{1}_{n_1} & 0 \\ 0 & 0 \end{pmatrix} \text{ and } \rho_{\tilde{\mathcal{B}}}(r_2) = E_2 = \begin{pmatrix} 0 & 0 \\ 0 & \mathbf{1}_{n_2} \end{pmatrix}.$$

Let $\mathcal{R}$ act on $\tilde{\mathcal{B}}^n$ via $\rho_{\tilde{\mathcal{B}}}$ and consider $(\varepsilon_1, \ldots, \varepsilon_n)$ its canonical basis. We define $\mathcal{L}$ as the $\mathcal{B}$-sublattice of $\tilde{\mathcal{B}}^n$ generated by $\rho_{\tilde{\mathcal{B}}}(r).\varepsilon_n$ when $r$ varies in $\mathcal{R}$. By construction, $\mathcal{L}$ is stable by the action of $\mathcal{R}$ and therefore is stable by the idempotents $E_i$. If we set $\mathcal{L}_i = E_i(\mathcal{L})$ then we have trivially $\mathcal{L} = \mathcal{L}_1 \oplus \mathcal{L}_2$.

**Lemma 1.1** $\mathcal{L}_2$ *is free of rank $n_2$ over $\mathcal{B}$.*

**Proof:** By definition of $\mathcal{L}_2$, one can see that $\mathcal{L}_2 \otimes_{\mathcal{B}} \kappa_{\mathcal{B}}$ is generated by

$$\begin{pmatrix} 0 & 0 \\ 0 & \bar{\rho}_2(r) \end{pmatrix} (\bar{\varepsilon}_n), \quad \forall r \in \mathcal{R}$$

and is therefore equal to $\kappa_{\mathcal{B}}.\bar{\varepsilon}_{n_1+1} \oplus \ldots \oplus \kappa_{\mathcal{B}}.\bar{\varepsilon}_n$. For each $i \in \{n_1+1, \ldots, n\}$, let $\varepsilon'_i \in \mathcal{L}_2$ such that $\bar{\varepsilon}'_i = \bar{\varepsilon}_i$. Then by Nakayama's lemma $\mathcal{L}_2 = \mathcal{B}.\varepsilon'_{n_1+1} + \ldots + \mathcal{B}.\varepsilon'_n$. Moreover this sum is direct since $\mathcal{L}_2 \otimes F_{\mathcal{B}}$ is of rank $n_2 = n - n_1$ over $F_{\mathcal{B}}$.∎

**Lemma 1.2** *The following map of $\mathcal{B}$-modules is surjective:*

$$\begin{aligned} \mathcal{R} &\to Hom_\mathcal{B}(\mathcal{L}_2, \mathcal{L}_1) \\ r &\mapsto Res_{\mathcal{L}_2}(E_1 \circ \rho_{\tilde{\mathcal{B}}}(r)) \end{aligned}$$

**Proof:** By Nakayama's lemma, we just have to prove that

$$\mathcal{R} \otimes_\mathcal{B} \kappa_\mathcal{B} \to Hom_{\kappa_\mathcal{B}}(\mathcal{L}_2 \otimes_\mathcal{B} \kappa_\mathcal{B}, \mathcal{L}_1 \otimes_\mathcal{B} \kappa_\mathcal{B}) \cong \bigoplus_{i=n_1+1}^{n} Hom_{\kappa_\mathcal{B}}(\kappa_\mathcal{B}.\bar{\varepsilon}_i, \mathcal{L}_1 \otimes_\mathcal{B} \kappa_\mathcal{B})$$

is surjective. By construction of $\mathcal{L}$, one sees immediately that the image projects surjectively on the last factor. If $x$ is any element of Since $\bar{\rho}_2$ is absolutely irreducible, for all $i \in \{n_1+1, \ldots, n\}$, there exist $s_i \in \mathcal{R}$ such that $\bar{\rho}_{\tilde{\mathcal{B}}}(r_2 s_i).\bar{\varepsilon}_n = \bar{\rho}_2(s_i).\bar{\varepsilon}_n = \bar{\varepsilon}_i$ and thus the image projects also surjectively on the other factors. Therefore considering $r \mapsto Res_{\mathcal{L}_2}(E_1 \circ \rho_{\tilde{\mathcal{B}}}(rs_i))$, one sees easily that the image contains the $i$-th factor for all $i$.∎

**Using the trace.** We are going to use the same method of the proof of the theorem of [31]. Thanks to the decomposition $\mathcal{L} = \mathcal{L}_1 \oplus \mathcal{L}_2$, we can see $\rho_\mathcal{L} = \rho_{\tilde{\mathcal{B}}}|_\mathcal{L}$ by blocks in the following manner:

$$\rho_\mathcal{L}(r) = \begin{pmatrix} A_r & B_r \\ C_r & D_r \end{pmatrix}$$

with

$$\begin{aligned} A_r &= Res_{\mathcal{L}_1}(E_1 \circ \rho_{\tilde{\mathcal{B}}}(r)) \in Hom_\mathcal{B}(\mathcal{L}_1, \mathcal{L}_1) \\ B_r &= Res_{\mathcal{L}_2}(E_1 \circ \rho_{\tilde{\mathcal{B}}}(r)) \in Hom_\mathcal{B}(\mathcal{L}_2, \mathcal{L}_1) \\ C_r &= Res_{\mathcal{L}_1}(E_2 \circ \rho_{\tilde{\mathcal{B}}}(r)) \in Hom_\mathcal{B}(\mathcal{L}_1, \mathcal{L}_2) \\ D_r &= Res_{\mathcal{L}_2}(E_2 \circ \rho_{\tilde{\mathcal{B}}}(r)) \in Hom_\mathcal{B}(\mathcal{L}_2, \mathcal{L}_2) \end{aligned}$$

Moreover, viewing these morphisms of $\mathcal{B}$-modules, as matrices with entries in $F_\mathcal{B}$, one can compute the traces of $A_r$, $D_r$ and $\rho_\mathcal{L}(r)$.

**Lemma 1.3** *For all $r \in \mathcal{R}$, $tr(A_r) \in \mathcal{B}$ and $tr(D_r) \in \mathcal{B}$. Moreover, for all $r \in \mathcal{R}$ we have:*

$$\begin{aligned} tr(A_r) &\equiv tr(\rho_1(r)) \bmod I \\ tr(D_r) &\equiv tr(\rho_2(r)) \bmod I \\ tr(C_r B_s) &\equiv tr(C_s B_r) \bmod I \end{aligned}$$

**Proof:** Let us prove it for $A_r$. We just have to remark that $tr(A_r) = tr(E_1 \rho_{\tilde{\mathcal{B}}}(r) E_1) = tr(\rho_{\tilde{\mathcal{B}}}(r_1 r r_1)) \in \mathcal{B}$. The congruences can be proven in the same manner of the proof of theorem of [31]( statements (Tr1) and (Tr2)).∎

We prove now by induction on $j \geq 1$ that $C_r \in Hom_\mathcal{B}(\mathcal{L}_1, (\mathcal{M}_\mathcal{B}^j + I)\mathcal{L}_2)$. It is obvious for $j = 0$. Assume that it is true for $j$ and let us prove it for $j + 1$.

**First Step:** Let $r \in Ker(\rho \otimes \kappa_\mathcal{B})$. Then $C_r \in Hom_\mathcal{B}(\mathcal{L}_1, (\mathcal{M}_\mathcal{B}^{j+1} + I)\mathcal{L}_2)$

**Proof:** Since $r \in Ker(\rho \otimes \kappa_\mathcal{B})$, we have $B_r(\mathcal{L}_2) \subset \mathcal{M}_\mathcal{B}.\mathcal{L}_1$, therefore for all $s \in \mathcal{R}$ by our induction hypothesis we have:

$$C_s B_r(\mathcal{L}_2) \subset C_s(\mathcal{M}_\mathcal{B}.\mathcal{L}_1) = \mathcal{M}_\mathcal{B}.C_s(\mathcal{L}_1) \subset (\mathcal{M}_\mathcal{B}^{j+1} + I).\mathcal{L}_2.$$

By lemma 1.3, we thus have $tr(C_r B_s) \in \mathcal{M}_\mathcal{B}^{j+1} + I$ for all $s \in \mathcal{R}$. Let $m \in \mathcal{L}_1$, we want to prove that $C_r.m \in (\mathcal{M}_\mathcal{B}^{j+1} + I)\mathcal{L}_2$. Let us write $C_r.m = \sum_{i=n_1+1}^n \alpha_i.\varepsilon_i'$ (with $\alpha_i \in \mathcal{B}$). For each $i \in \{n_1 + 1, \ldots, n\}$, by lemma 1.2, there exists $t_i \in \mathcal{R}$ such that $B_{t_i}(\varepsilon_i') = m$ and $B_{t_i}(\varepsilon_{i'}') = 0$ if $i' \neq i$. Then $\alpha_i = tr(C_r B_{t_i}) \in \mathcal{M}_\mathcal{B}^{j+1} + I$.■

**Second Step:** We proceed as in the second part of the proof of the theorem of [31] in order to prove that $C_r \in Hom_\mathcal{B}(\mathcal{L}_1, (\mathcal{M}_\mathcal{B}^{j+1} + I)\mathcal{L}_2)$ for all $r \in \mathcal{R}$. Consider $Im\rho \otimes \kappa_\mathcal{B} \subset Hom(\mathcal{L} \otimes \kappa_\mathcal{B}, \mathcal{L} \otimes \kappa_\mathcal{B})$. We denote by $\bar{A}_r, \bar{B}_r, \ldots$ the projections of $\bar{\rho}(r)$ on $Hom(\mathcal{L}_1 \otimes \kappa_B r, \mathcal{L}_1 \otimes \kappa_B r)$, $Hom(\mathcal{L}_2 \otimes \kappa_B r, \mathcal{L}_1 \otimes \kappa_B r)$,...These maps induce (by using the projectors $E_1$ and $E_2$) a decomposition $Im\rho \otimes \kappa_\mathcal{B} = (Im\rho \otimes \kappa_\mathcal{B})_{11} \oplus (Im\rho \otimes \kappa_\mathcal{B})_{12} \oplus (Im\rho \otimes \kappa_\mathcal{B})_{21} \oplus (Im\rho \otimes \kappa_\mathcal{B})_{22}$ and we will denote an element of this image by a matrix $\begin{pmatrix} \bar{A} & \bar{B} \\ \bar{C} & \bar{D} \end{pmatrix}$ with $A \in (Im\rho \otimes \kappa_\mathcal{B})_{11}$, $B \in (Im\rho \otimes \kappa_\mathcal{B})_{12}$...From the first step, we have a map $\Phi$ from $Im\rho \otimes \kappa_\mathcal{B}$ to $Hom_\mathcal{B}(\mathcal{L}_1, (\mathcal{M}_\mathcal{B}^j + I)\mathcal{L}_2/, (\mathcal{M}_\mathcal{B}^{j+1} + I)\mathcal{L}_2)$ induced by $r \mapsto C_r$. Moreover by the relation $C_{rs} = C_r A_s + D_r C_s$, we see that

$$\Phi(\begin{pmatrix} \bar{A} & \bar{B} \\ \bar{C} & \bar{D} \end{pmatrix} \begin{pmatrix} \bar{A}' & \bar{B}' \\ \bar{C}' & \bar{D}' \end{pmatrix}) = \phi(\begin{pmatrix} \bar{A} & \bar{B} \\ \bar{C} & \bar{D} \end{pmatrix}) \circ \bar{A}' + \bar{D} \circ \Phi(\begin{pmatrix} \bar{A}' & \bar{B}' \\ \bar{C}' & \bar{D}' \end{pmatrix})$$

By using the fact that $\Phi(\bar{\rho}(r_1)) = \Phi(\bar{\rho}(r_2)) = 0$ and the above relation, it is then easy to see that $\Phi = 0$ and the induction is proved for $j + 1$.■

The consequence of that fact is that $\mathcal{L}_1 \otimes \mathcal{B}/I$ is stable by the action of $\mathcal{R}$ and the action on the quotient $(\mathcal{L} \otimes \mathcal{B}/I)/(\mathcal{L}_1 \otimes \mathcal{B}/I)$ is isomorphic to $\rho_2 \otimes \mathcal{B}/I$ by Carayol's Theorem and lemma 1.3. In order to conclude the proof of our theorem, let us examine the action on $\mathcal{L}_1 \otimes \mathcal{B}/I$.

**Lemma 1.4** *The $\mathcal{B}$-modules $\mathcal{L}_1(\alpha_i) = E(\alpha_i)\mathcal{L}_1 = Ker(\rho_{\tilde{\mathcal{B}}}(r_0) - \alpha_i Id) \cap \mathcal{L}_1$ for $i \in \{1, \ldots, n_1\}$ are mutually isomorphic.*

**Proof:** Recall that for all $i$, $\rho_{\tilde{\mathcal{B}}}(s_i)$ is the matrix having zero everywhere but a 1 at the entry (i,i). For $i, i' \in \{1, \ldots, n_1\}$, there exist by irreducibility of $\bar{\rho}_1$ an element $\tau_{i,i'} \in \mathcal{R}$ such that $\bar{\rho}_{\tilde{\mathcal{B}}}(\tau_{i,i}).\bar{\varepsilon}_i = \bar{\varepsilon}_{i'}$. Then taking $\rho_{\tilde{\mathcal{B}}}(s_{i'}\tau_{i',i}s_i)$, we get $\phi_{i,i'}$ a morphism $\mathcal{L}_1(\alpha_i)$ into $\mathcal{L}_1(\alpha_{i'})$. By our choice of the $\tau_{i',i}$'s, $\phi_{i,i'} \circ \phi_{i',i}$ is the identity on $Ker(\rho_{\tilde{\mathcal{B}}}(r_0) - \alpha_i Id) \otimes \kappa_{\mathcal{B}}$ and thus is bijective on $Ker(\rho_{\tilde{\mathcal{B}}}(r_0) - \alpha_i Id)$; since it respects $\mathcal{L}_1(\alpha_i)$ it is an automorphism of $\mathcal{L}_1(\alpha_i)$. Therefore the $\phi_{i,i'}$ are isomorphisms.∎

Let us set $\mathcal{T} = \mathcal{L}_1(\alpha_1)$ and fix an isomorphism $\mathcal{L}_1 \otimes \mathcal{B}/I \cong (\mathcal{T}/I\mathcal{T})^{n_1}$, and consider the non-commutative artinian algebra $\mathcal{E} = End_{\mathcal{B}/I}(\mathcal{T}/I\mathcal{T})$ denoting by $\theta$ the canonical homomorphism of algebra $\mathcal{B}/I \to \mathcal{E}$. Then the action of $\mathcal{R}$ on $\mathcal{L}_1 \otimes \mathcal{B}/I$ gives us a representation $\rho'_1$ in $M_{n_1}(\mathcal{E})$ such that $tr(\rho'_1(r)) \in \theta(\mathcal{B}/I)$ is defined for all $r \in \mathcal{R}$ and such that $tr(\rho'_1(r)) = \theta(tr(\rho_1(r)))$ for all $r \in \mathcal{R}$. Then by a simple generalization of Carayol's theorem (cf. proof 1.1.2 of [3]), we get that

$$\mathcal{L}_1 \otimes \mathcal{B}/I \cong \rho_1 \otimes \mathcal{T}/I\mathcal{T}.$$

This point finishes the proof of theorem 1.1.∎

## 1.2 Applications to representations in $GSp_4$.

Let

$$GSp_4 = \{g \in GL_4; {}^t g J g = \nu(g) J\} \text{ where } J = \begin{pmatrix} 0_2 & -1_2 \\ 1_2 & 0_2 \end{pmatrix}.$$

Let $G$ be a group and $\rho$ be an absolutely irreducible representation of $G$ in $GSp_4(F_{\mathcal{B}})$ and let $\nu_\rho$ be the composite of $\rho$ and the multiplier character of $GSp_4$. Assume that there exist a representation $\rho_0 : G \to GL_2(\mathcal{B})$, a character $\nu_0 : G \to \mathcal{B}^\times$ and $I \subset \mathcal{B}$ an ideal of $\mathcal{B}$ such that:

(i) The coefficients of the characteristic polynomial of $\rho$ belongs to $\mathcal{B}$.

(ii) The characteristic polynomials of $\rho$ and $\rho_0 \oplus {}^t\rho_0^{-1} \otimes \nu_0$ are congruent modulo $I$.

(iii) $\bar{\rho}_0$ is absolutely irreducible

(iv) $\bar{\nu}_0 \neq det(\bar{\rho}_0)$

Then by theorem 1.1, there exists a $\rho(G)$-stable $\mathcal{B}$-lattice $\mathcal{L} = \mathcal{L}_1 \oplus \mathcal{L}_2$ in $F_{\mathcal{B}}^4$ and a finitely generated $\mathcal{B}$-submodule $\mathcal{T}$ of $F_{\mathcal{B}}$ such that $\mathcal{L}_1 \cong \mathcal{B}^2 \otimes \mathcal{T}$ and $\mathcal{L}_1 \cong \mathcal{B}^2$ such that the exact sequence of $\mathcal{B}$-module

$$0 \to \mathcal{L}_1 \to \mathcal{L} \to \mathcal{L}_2 \to 0$$

gives us modulo $I$ the following exact sequence of $\mathcal{B}[G]$-module.

$$0 \longrightarrow V(\rho_0) \otimes \mathcal{T}/I.\mathcal{T} \longrightarrow \mathcal{L} \otimes \mathcal{B}/I \xrightarrow{s} V({}^t\rho_0^{-1} \otimes \nu_0) \otimes \mathcal{B}/I \longrightarrow 0$$

where $V(\rho_0) = \mathcal{B}^2$ with $G$-action given by $\rho_0$. Here $s$ is not $G$-linear; more precisely we know that $\mathcal{L}$ has no quotient isomorphic to $\bar{\rho}_0$.

Let $\phi$ be a "$\nu_\rho$-invariant" skew-symmetric $F_\mathcal{B}$-valued bilinear form on $F_\mathcal{B}^4$ (it exists since $\rho$ takes values in $GSp_4(F_\mathcal{B})$) and let us consider the restriction $\phi_\mathcal{L}$ of $\phi$ to $\mathcal{L} \otimes_\mathcal{B} \mathcal{L}$ and set $\mathcal{S} = \phi_\mathcal{L}(\mathcal{L} \otimes_\mathcal{B} \mathcal{L})$. Then we have:

• $\phi_\mathcal{L}(\mathcal{L}_1 \otimes \mathcal{L}_1) \subset I.\mathcal{S}$: Indeed, let $v, w \in \mathcal{L}_1$ and $\sigma_0 \in G$ such that $det(\rho_0(\sigma_0)) \not\equiv \nu_0(\sigma_0); mod\ \mathcal{M}_\mathcal{B}$. Then $\nu_0(\sigma_0)\phi_\mathcal{L}(v,w) \equiv \nu_\rho(\sigma_0)\phi_\mathcal{L}(v,w) = \phi_\mathcal{L}(\sigma_0.v, \sigma_0.w) \equiv det(\bar{\rho}_0(\sigma_0))\phi_\mathcal{L}(v,w)$ mod $I.\mathcal{S}$, thus $\phi_\mathcal{L}(\mathcal{L}_1 \otimes_\mathcal{B} \mathcal{L}_1) \subset I.\mathcal{S}$. By similar arguments, we can prove that:

• $\phi_\mathcal{L}(\mathcal{L}_2 \otimes \mathcal{L}_2) \subset I.\mathcal{S} + \phi_\mathcal{L}(\mathcal{L}_1 \otimes \mathcal{L}_2)$.

• Study of $\phi_\mathcal{L}(\mathcal{L}_1 \otimes \mathcal{L}_2)$: Let $\sigma_1 \in G$ such that $\rho_0(\sigma_1)$ has two eigenvalues $\alpha_1$ and $\alpha_2$ in $\mathcal{B}$ which are distinct modulo $\mathcal{M}_\mathcal{B}$. By the above exact sequence, there exist $e_1, e_2 \in \mathcal{B}^2$ such that $\mathcal{L}_1 = \mathcal{B}.e_1 \otimes \mathcal{T} \oplus \mathcal{B}.e_2 \otimes \mathcal{T}$ and such that $\sigma_1.e_i \equiv \alpha_i.e_i$ mod. $I.\mathcal{L}_1$. There exists also a basis $(e'_1, e'_2)$ of $\mathcal{L}_2$ such that $\sigma_1.e'_i \equiv \nu_0(\sigma_1)\alpha_i^{-1}.e'_i$ mod. $I.\mathcal{L}_2 + \mathcal{L}_1$. Then for all $t \in \mathcal{T}$,

$$\nu_0(\sigma_1)\phi_\mathcal{L}(e_1 \otimes t, e'_2) \equiv \nu_\rho(\sigma_1)\phi_\mathcal{L}(e_1 \otimes t, e'_2) = \\ \phi_\mathcal{L}(\sigma_1.e_1 \otimes t, \sigma_1.e'_2) \equiv \alpha_1\nu_0(\sigma_1)\alpha_2^{-1}\phi_\mathcal{L}(e_1 \otimes t, e'_2)$$

and thus $\phi_\mathcal{L}(e_1 \otimes t, e'_2) \in I.\mathcal{S}$ since $\alpha_1\alpha_2^{-1} \not\equiv 1$ mod. $\mathcal{M}_\mathcal{B}$ (and also $\phi_\mathcal{L}(e_2 \otimes t, e'_1) \in I.\mathcal{S}$). Let now be $\sigma_2 \in G$ such that

$$Mat_{(e_1\ mod.\ I.\mathcal{L}, e_2\ mod.\ I.\mathcal{L})}(\sigma_2) = \begin{pmatrix} a & b \\ c & d \end{pmatrix}$$

with $bc \notin \mathcal{M}_\mathcal{B}$ (it is possible since $\bar{\rho}_0$ is absolutely irreducible). After multiplying (if necessary) $e'_1$ and $e'_2$ by suitable elements of $\mathcal{B}$, we have also:

$$Mat_{(e'_1\ mod.\ I.\mathcal{L}, e'_2\ mod.\ I.\mathcal{L})}(\sigma_2) = \frac{\nu_0(\sigma_2)}{ad-bc}\begin{pmatrix} d & -c \\ -b & a \end{pmatrix}.$$

then, for all $t \in \mathcal{T}$

$$\nu_0(\sigma_2)\phi_\mathcal{L}(e_1 \otimes t, e'_1) \equiv \\ \phi_\mathcal{L}(a.e_1 \otimes t + c.e_2 \otimes t, \tfrac{\bar{\nu}_0(\sigma_2)}{ad-bc}(d.e'_1 - b.e'_2))\ mod.\ I.\mathcal{S}$$

and thus $\phi_\mathcal{L}(e_1 \otimes t, e'_1) \equiv \phi_\mathcal{L}(e_2 \otimes t, e'_2)$ mod. $I.\mathcal{S}$ since $bc \notin \mathcal{M}_\mathcal{B}$. Therefore we have proven that $\mathcal{S} = \phi_\mathcal{L}((e_1 \otimes \mathcal{T}), e'_1) + I.\mathcal{S}$. By Nakayama's lemma, we thus have $\mathcal{S} = \phi_\mathcal{L}((e_1 \otimes \mathcal{T}), e'_1) \cong \mathcal{T}$ by the isomorphism $t \mapsto \phi_\mathcal{L}(e_1 \otimes t, e'_1)$. Considering $\phi_\mathcal{L}$ modulo $I$, we get from the above discussion the following:

**Proposition 1.1** *There exists a skew-symmetric $\mathcal{T}/I\mathcal{T}$-valued bilinear form $\bar{\phi}$ on $\mathcal{L}/I\mathcal{L}$ such that*

$$\phi_\mathcal{L}(g.v, g.w) = \nu_0(g)\phi_\mathcal{L}(v,w)$$

Moreover $V(\rho_0) \otimes \mathcal{T}/I.\mathcal{T}$ and $s(V({}^t\rho_0^{-1} \otimes \nu_0) \otimes \mathcal{B}/I)$ are isotropic for $\bar{\phi}$ and

$$\bar{\phi}(\begin{pmatrix} t_1 \\ t_2 \end{pmatrix}, \begin{pmatrix} s(b_1) \\ s(b_2) \end{pmatrix})) = b_1.t_1 + b_2.t_2$$

for all $\begin{pmatrix} t_1 \\ t_2 \end{pmatrix} \in V(\rho_0) \otimes \mathcal{T}/I$ and $\begin{pmatrix} b_1 \\ b_2 \end{pmatrix} \in V({}^t\rho_0^{-1} \otimes \nu_0) \otimes \mathcal{B}/I$.

## 2 The Eisenstein-Klingen Ideal

We begin by giving some review on the ordinary universal Hecke algebra for $GSp_4/\mathbf{Q}$ studied in [27],[28].

### 2.1 Notations and definitions

Let $\mathbf{A}$ (resp. $\mathbf{A}_f$) be the ring of adeles (resp. of finite adeles) of $\mathbf{Q}$.
Let

$$G = GSp_4 = \{g \in GL_4; {}^tgJg = \nu(g)J\} \text{ where } J = \begin{pmatrix} 0_2 & -1_2 \\ 1_2 & 0_2 \end{pmatrix}$$

and $G' = Sp_4$ its derived subgroup.

Let us introduce the standard Borel subgroup of $GSp_4$:

$$B = \{\begin{pmatrix} \times & \times & \times & \times \\ 0 & \times & \times & \times \\ 0 & 0 & \times & 0 \\ 0 & 0 & \times & \times \end{pmatrix} \in GSp_4\}$$

Let $U$ be the unipotent radical of $B$ and $T$ be the maximal torus contained in $B$. For any elements $a, b, c$ in any ring $A$, we put

$$[a, b; c] = Diag(a, b, ca^{-1}, cb^{-1}) \in T(A)$$

For $m, n \in \mathbf{Z}$ we define the algebraic character $\chi_{(m,n)}$ of $T$ by:

$$\chi_{(m,n)}([a, b; c]) = a^m b^n.$$

It is called dominant (resp. Klingen-regular dominant, Siegel-regular dominant or regular dominant) if $m \geq n \geq 0$ (resp. $m \geq n > 0$, $m > n \geq 0$ or $m > n > 0$).

**Definition 2.1** *An $O_K$-valued character $\chi$ of $T(\mathbf{Z}_p)$ will be called arithmetic of level $p^r$ of weight $(m, n)$, if it is equal to $\chi_{(m,n)}$ on the kernel of $T(\mathbf{Z}_p) \to T(\mathbf{Z}/p^r\mathbf{Z})$. That means that $\chi = \chi_{(m,n)} \times \epsilon$ where $\epsilon$ is a finite character of level $p^r.\chi_{(m,n)}$ is called the algebraic part of $\chi$ and is denoted by $\chi^{alg}$.*

## 2.2 p-Ordinary cohomology of the Siegel threefold

Let $K_\infty$ be the stabilizer of the map $h : \mathbf{C}^\times \longrightarrow G(\mathbf{R})$ given by

$$h(x+iy) = \begin{pmatrix} x 1_2 & y 1_2 \\ -y 1_2 & x 1_2 \end{pmatrix}.$$

and $K'_\infty = K_\infty \cap G'(\mathbf{R})$. For any compact open subgroup $V$ of $G'(\mathbf{A}_f)$, we consider the connected Siegel threefold:

$$X(V) = G'(\mathbf{Q}) \backslash G'(\mathbf{A}) / K'_\infty . V$$

For $r \geq 1$, we set

$$U_1(Np^r) = \{\gamma \in G'(\hat{\mathbf{Z}}) \text{ such that } \gamma \bmod p^r \in U_B(\mathbf{Z}/p^r\mathbf{Z})\}$$
$$U_0(Np^r) = \{\gamma \in G'(\hat{\mathbf{Z}}) \text{ such that } \gamma \bmod p^r \in B(\mathbf{Z}/p^r\mathbf{Z})\}$$

For any module $L$ over which $G'(\mathbf{Z}_p)$ acts linearly, one considers its associated local system as the shief of locally constant sections of the following cover:

$$G'(\mathbf{Q}) \backslash G(\mathbf{A}) \times L / K'_\infty . V \to G'(\mathbf{Q}) \backslash G(\mathbf{A}) / K'_\infty . V = X(V)$$

where $\gamma.(g,m).(k_\infty, v) = (\gamma.g.(k_\infty, v), v^{-1}.m)$.

We describe now the modules we are interested in. We fix $K$ a finite extension of $\mathbf{Q}_p$ and let $O_K$ (resp. $\varpi$ and $\kappa$) be the ring of integers in $K$ (resp. a uniformizing element and the residue field of $O_K$). The algebraic representation of $G_{/\mathbf{Q}_p}$ of highest weight $\chi_{(m,n)}$ can be describe as:

$$L_{(m,n)}(K) = \{f : G(\mathbf{Q}_p) \to K; f(tg) = \chi_{(m,n)}(t^{-1})f(g)\},$$

the action being given by $(\gamma.f)(g) = f(\gamma^{-1}g)$.

Let $I$ be the Iwahori subgroup of $G(\mathbf{Q}_p)$ (i.e. the elements $g$ of $G(\mathbf{Z}_p)$ such that $g \bmod p$ belongs to $B(\mathbf{Z}/p\mathbf{Z})$).

For any $O_K$-valued arithmetic character $\chi$ of level $p^r$ and weight $(m,n)$, we consider the following integral structure of the above representation space $L_{(m,n)}(K)$:

$$L_\chi(O_K) = \{f \in L_{(m,n)}(K); f(I) \subset O_K\}$$

where the subscript $\chi$ means that the natural induced action of $I$ is twisted by the finite character $\epsilon$ for which $\chi = \chi_{(m,n)} \times \epsilon$.

Then we will denote by $\mathcal{L}_\chi$ the local system associated to the representation $L_\chi = L_\chi(O_K) \otimes_{O_K} K/O_K$..

For $U_r = U_0(Np^r)$ or $U_1(Np^r)$, we are interested in the cohomology $H^3_!(X(U_r); \mathcal{L}_\chi)$ where $H_!$ means the image of compact support cohomology in the total cohomology. More Precisely, we can study only the ordinary part of that cohomology. Let us now define this notion in that context. We consider the double classes

$$T_{1,p} = U_r[1,1;p]U_r$$
$$T_{2,p} = U_r[1,p;p^2]U_r$$

Since $L_\chi$ is a cofinite $O_K$-module over which $[1,1;p]$ and $[1,p;p^2]$ act, the above double classes act on the cohomology of $\mathcal{L}_\chi$ and we can consider the ordinary idempotent (cf. [27] and [28]) acting on the cohomology we are interested in

$$e_p = \lim_{n \to \infty} (T_{1,p}T_{2,p})^{n!}$$

and set

$$\mathcal{V}_\chi(Np^r) = e_p.H^3_!(X(U_0(Np^r));\mathcal{L}_\chi)$$

and

$$\mathcal{V} = \mathcal{V}(N) = \varinjlim_r e_p.H^3_!(X(U_1(Np^r));\widetilde{K/O_K})$$

We let act $T'(\mathbf{Z}_p)$ on these groups by the double classes $<t> = [U_r t U_r]$ for all $t \in T'(\mathbf{Z}_p)$. For any arithmetic character $\chi$, we set

$$\mathcal{V}[\chi] = \{v \in \mathcal{V}; t.v = \chi(t)v \; \forall t \in T'(\mathbf{Z}_p)\}$$

By the results of [28], we have:

**Theorem 2.1** *Let $\chi$ be an arithmetic $O_K$-valued characters of $T(\mathbf{Z}_p)$ of level $p^r$. If $\chi$ is dominant and regular, then the canonical following map is an isogeny.*

$$\mathcal{V}_\chi(Np^r) \to \mathcal{V}[\chi]$$

The decomposition $\mathbf{Z}_p^\times = \mu_{p-1} \times (1+p\mathbf{Z}_p)$ yields a decomposition $T'(\mathbf{Z}_p) = \Delta \times W$ where $\Delta$ is the torsion part of $T'(\mathbf{Z}_p)$ and $W$ is rank 2 over $\mathbf{Z}_p$. Then we set $\mathbf{\Lambda} = O_K[[W]]$ and $\underline{\mathbf{\Lambda}} = O_K[[T'(\mathbf{Z}_p)]] = \mathbf{\Lambda}[\Delta]$. If we fix $u$ a topological generator of $(1+p\mathbf{Z}_p) \subset \mathbf{Z}_p^\times$ and set $T_1 = [u,1;1]-1$ and $T_2 = [1,u;1]-1$, we can identify $\mathbf{\Lambda}$ to the 2-variable power series ring $O_K[[T_1,T_2]]$.

**Corollary 2.1** $\mathbf{V}(N) = \mathcal{V}(N)^* = Hom_{\mathbf{Z}_p}(\mathcal{V}(N);\mathbf{Q}_p/\mathbf{Z}_p)$ *is of finite type over $\mathbf{\Lambda}$.*

## 2.3 The p-ordinary universal Hecke algebra

Let $\omega : \mathbf{Z}_p^\times \to \mathbf{Z}_p^\times$ be the Teichmüller character. We set for any $z \in \mathbf{Z}_p^\times$ and any variable $S$:

$$<z>_S = (1+S)^{\frac{\log_p(z\omega(z)^{-1})}{\log_p(u)}} \in \mathbf{Z}_p[[S]].$$

For any character $\xi$ of a group $\Gamma$ in $\mathbf{Z}_p^\times$, we denote by $<\xi>_S$ the composite of $<->_S \circ \xi$.

For any arithmetic character $\chi$, we denote by $P_\chi$ (resp. $\underline{P}_\chi$) the kernel of the canonical homomorphism from $\mathbf{\Lambda}$ (resp. from $\underline{\mathbf{\Lambda}}$) into $O_K$ induced by $\chi$, and $\omega_{m,n}$ the restriction of $\chi_{(m,n)}$ to $\Delta$ ( note that it depends only on the class of $(m,n)$ in $(\mathbf{Z}/(p-1)\mathbf{Z})^2$).

**Definition 2.2** *(i) An irreducible closed subscheme $\mathrm{Spec}(A)$ of $\mathrm{Spec}(\mathbf{\Lambda})$ is called arithmetic if and only if $\mathrm{Spec}(A)(\bar{\mathbf{Q}}_p) \subset \mathrm{Spec}(\mathbf{\Lambda})(\bar{\mathbf{Q}}_p)$ contains at least (in $\mathrm{Spec}(A)$) one arithmetic point. It is called ?-regular if $\mathrm{Spec}(A)(\bar{\mathbf{Q}}_p)$ contains at least one arithmetic ?-regular point.*

*(ii) Any finite and flat irreducible scheme $\mathrm{Spec}(\mathbf{J})$ over an arithmetic (resp. arithmetic ?-regular) irreducible closed subscheme of $\mathrm{Spec}(\mathbf{\Lambda})$ will be called arithmetic (resp. arithmetic ?-regular).*

For any $\ell \not| Np$, we consider the double classes

$$T_\ell = U_r[1,1;\ell]U_r$$
$$R_\ell = U_r[1,\ell;\ell^2]U_r$$
$$S_\ell = U_r[\ell,\ell;\ell^2]U_r$$

As it is well known, these double classes act on $\mathcal{V}_\chi(Np^r)$ and $\mathcal{V}(N)$. We denote by $h_\chi(Np^r)$ ( respectively $\mathbf{h} = \mathbf{h}(N)$) the $O_K$-algebra of $End_{O_K}(\mathcal{V}_\chi(Np^r))$ ( respectively $End_{O_K}(\mathcal{V}(N))$) generated by the image of $T_\ell$, $R_\ell$, $S_\ell$ and $<t>$ for all $\ell \not| N$ and $t \in T'(\mathbf{Z}_p)$. We set $h_\chi(Np^r;K) = h_\chi(Np^r) \otimes K$. We will consider also $R(N;O_K)$ the abstract Hecke algebra generated over $O_K$ by the variables $T_\ell$, $R_\ell$, $S_\ell$ and $<t>$ for all $\ell \not| N$ and $t \in T'(\mathbf{Z}_p)$.

**Corollary 2.2** *Let $\chi$ be an arithmetic $O_K$-valued characters of $T'(\mathbf{Z}_p)$ of level $p^r$ and weight $(m,n)$. If $\chi$ is dominant and sufficiently regular, and let $P = \underline{P}_\chi$,*

$$\mathbf{h} \otimes (\underline{\mathbf{\Lambda}})_P/P \to h_\chi(Np^r;K)$$

*is surjective and its kernel is contained in the radical of $\mathbf{h} \otimes (\underline{\mathbf{\Lambda}})_P/P$.*

It is a classical fact that if $\lambda$ is a character of the Hecke algebra $h_\chi(Np^r; K)$ in $K$ there exists a cohomological irreducible cuspidal representation of $GSp_4/\mathbf{Q}$ of cohomological weight $\chi^{alg}$ of level $Np^r$ whose Langlands' parameters at $\ell \nmid Np$ are given by the roots of the polynomial $\lambda(\mathcal{Q}_\ell)$ where

$$\mathcal{Q}_\ell = X^4 - T_\ell X^3 + \ell(R_\ell + (1+\ell^2)S_\ell)X^2 - \ell^3 T_\ell S_\ell X + \ell^6 S_\ell^2$$

Such a representation is called associated to $\lambda$. Conversely, such a representation $\pi$ give us a character of the Hecke algebra $R(N, O_K)$ noted $\lambda_\pi$ such that $\pi$ is associated to $\lambda_\pi$. In that case, $\pi$ will be called $p$-ordinary if $\lambda_\pi$ factorizes through $h_\chi(Np^r; K)$ for some $r$, $N$ and $\chi$.

Let $\mathbf{J}$ be a finite and flat extension of $\mathbf{\Lambda}$. If $\mathbf{J}$ is an irreducible component of the Hecke algebra $\mathbf{h}(N)$, we denote by $\lambda_\mathbf{J}$ the corresponding character of $\mathbf{h}(N)$. For $(a,b) \in (\mathbf{Z}/(p-1)\mathbf{Z})^2$, we say that $\mathbf{J}$ is of nebentypus $\omega_{a,b}$ if $\mathbf{J}$ is an irreducible component of $\mathbf{h}[\omega_{a,b}]$, the $\omega_{a,b}$ part under the action of $\Delta$ on $\mathbf{h}$. Therefore from corollary 2.2, we deduce easily:

**Corollary 2.3** *Let $\mathbf{J}$ be as above. For any prime $P$ of $\mathbf{J}$ above $P_\chi$ for dominant and sufficiently regular $\chi$ of level $p^r$, $\lambda_\mathbf{J}$ mod. $P$ factorizes through a character of $h_\chi(Np^r; K)$. i.e. There exists a character $\lambda_\chi$ such that the following diagram commutes:*

$$\begin{array}{ccc}
\mathbf{h}(N) & \xrightarrow{\lambda_\mathbf{J}} & \mathbf{J} \\
\downarrow & & \downarrow \\
\mathbf{h}(N) \otimes (\mathbf{\Lambda})_{P_\chi}/P_\chi & \longrightarrow & \mathbf{J}_P/P \\
\downarrow & & \| \\
h_{\chi \omega_{m,n}^{-1}}(Np^r; K) & \xrightarrow{\lambda_\chi} & \mathbf{J}_P/P
\end{array}$$

**Corollary 2.4** *Let $\mathrm{Spec}(A)$ be an irreducible arithmetic closed subscheme of $Spec(\mathbf{\Lambda})$. Let $F_A$ be the field of fractions of $A$, then $\mathbf{h} \otimes F_A$ is semisimple.*

**Proof:** We have to prove that any nilpotent element of $\mathbf{h} \otimes A$ is torsion. Let $\mathbf{T}$ be such an element and $P_\chi$ an arithmetic points of $\mathrm{Spec}A$. Suppose $\mathbf{T}$ is not torsion. Then we can assume that $\mathbf{T}.\mathcal{V} \otimes A \notin P_\chi.\mathcal{V} \otimes A$. However the Hecke operators act semi-simply on $V_\chi(Np^r) \otimes K$; therefore $\mathbf{T}$ should act trivially on that module. But this contredicts our assumption on $\mathbf{T}$ since $(V(N) \otimes A/P_\chi) \otimes K \cong V_\chi(Np^r) \otimes K$. Therefore $T$ is torsion in $\mathbf{h} \otimes A$.∎

## 2.4 The Eisenstein-Klingen character

First, let us recall the definition of the so-called Klingen parabolic subgroup of $GSp_4$.

$$Q = \{ \begin{pmatrix} \times & \times & \times & \times \\ 0 & \times & \times & \times \\ 0 & 0 & \times & 0 \\ 0 & \times & \times & \times \end{pmatrix} \in GSp_4 \}$$

The Levi subgroup of $Q$ is isomorphic to $\mathbf{G}_m \times GL_2$ by

$$(e, \begin{pmatrix} a & b \\ c & d \end{pmatrix}) \mapsto \begin{pmatrix} e & 0 & 0 & 0 \\ 0 & a & 0 & b \\ 0 & 0 & e^{-1}(ad-bc) & 0 \\ 0 & c & 0 & d \end{pmatrix}$$

The proof of the following lemma is an easy computation on Langlands' parameters associated to unramified representations.

**Lemma 2.1** *Let $E$ be a non archimedean local fields. Let $\sigma_E$ (resp. $\chi_E$) be an unramified representation of $GL_{2/E}$ (resp. an unramified character of $E^\times$) and $\pi_E$ the Langlands quotient of $Ind_{Q(E}^{G(E)} \sigma_E \otimes \chi_E$. Then we have*

$$L(s, \pi_E) = L(s, \sigma_E) \times L(s-2, \sigma_E \otimes \chi_E)$$

Let now $\sigma$ be a cuspidal representation of $GL_2/\mathbf{Q}$ of level $N$ whose archimedean component is the discrete series representation $\pi(k-1, 2-k)$ of $GL_2(\mathbf{R})$ corresponding to the Weil parametrization:

$$W_{\mathbf{R}} = Gal(\mathbf{C}/\mathbf{R}) \times \mathbf{C}^\times \rightarrow GL_2(\mathbf{C})$$

$$z \in \mathbf{C} \mapsto \begin{pmatrix} z^{1-k}|z|^{1/2} & 0 \\ 0 & \bar{z}^{1-k}|z|^{1/2} \end{pmatrix}$$

$$c \mapsto \begin{pmatrix} 0 & 1 \\ -1 & 0 \end{pmatrix}.$$

Let $\epsilon$ be a Dirichlet character modulo $N$ and $t \in \mathbf{Z}$. We consider the parabolic induction:

$$Ind_{Q(\mathbf{A})}^{G(\mathbf{A})} \sigma \otimes \epsilon |.|_{\mathbf{A}}^{-t}.$$

If $k - 3 \geq t$ by result of [1], the archimedean component of this induction has a quotient isomorphic to the discrete series representation of cohomological parameters $(k-3, t; 3-t-k)$. Therefore we can take a convenient $\phi \in Ind_{Q(\mathbf{A})}^{G(\mathbf{A})} \sigma \otimes \epsilon |.|_{\mathbf{A}}^{-t}$ and form Eisenstein series:

$$E_\phi(g) = \sum_{\gamma \in Q(\mathbf{Q}) \backslash G(\mathbf{Q})} \phi(\gamma g)$$

If $E_\phi$ converges (especially when the weight is sufficiently regular), it should produce cohomology class for the local system $L_{(k-3,t)}$. Moreover, by Lemma 2.1 the corresponding character of the Hecke algebra $\lambda_{\sigma,\epsilon|.|^{-t}}$ is such that

$$\lambda_{\sigma,\epsilon|.|^{-t}}(\mathcal{Q}_\ell(\ell^{-s})) = L(s,\sigma_\ell) \times L(s-2,\sigma_\ell \otimes \chi_\ell)$$

We can now attach an Eisenstein ideal to a couple $(\mathcal{F},\eta)$ where $\eta$ is a Dirichlet character of level $Np$ and $\mathcal{F}$ a Hida family of cuspidal ordinary elliptic cusp forms. Let us precise the data of that family. Let $\mathbf{I}$ be a finite and flat extension of $O_K[[T]]$. We take $\mathcal{F} = \sum_n a(n;F)q^n$ an $\mathbf{I}$-adic ordinary elliptic cusp form of tame level N and nebentypus $\varepsilon_\mathcal{F}$ (cf. [13] for more details). Then for such a couple $(\mathcal{F},\eta)$, we consider the character $\lambda_{\mathcal{F},\eta}$ of $R(N;O_K)$ in $\mathbf{I}[[S]]$ such that:

$$\lambda_{\mathcal{F},\eta}(\mathcal{Q}_\ell(X)) = (X^2 - a(\ell;F)X + \ell^{-1}\varepsilon_\mathcal{F}(\ell)<\ell>_T$$
$$\times (X^2 - a(\ell;\mathcal{F})\eta(\ell)<\ell>_S X + \ell^{-1}\varepsilon_\mathcal{F}(\ell)<\ell>_T \eta(\ell)^2 <\ell^2>_S)$$

and $\lambda_{\mathcal{F},\eta}(<[x,y;1]>) = <x>_{(1+T)u^{-3}-1} <y>_{(1+S)u^{-1}-1} \varepsilon_\mathcal{F}\omega^{-3}(x)\eta\omega^{-1}(y)$. We look at $\mathbf{I}[[S]]$ as a finite flat extension of $\Lambda$ via the map

$$\begin{aligned} \Lambda = O[[T_1,T_2]] &\xrightarrow{\iota} \mathbf{I}[[S]] \\ T_1 &\mapsto (1+T)u^{-3}-1 \\ T_2 &\mapsto (1+S)u^{-1}-1 \end{aligned}$$

Then we define $I_{Eis(\mathcal{F},\eta)}$ as the ideal of the Hecke algebra $\mathbf{h}(N) \otimes_\iota \mathbf{I}[[S]]$ generated by $1 \otimes \lambda_{\mathcal{F},\eta}(\mathbf{T}) - \mathbf{T} \otimes 1$ for all $\mathbf{T} \in R(N;O_K)$. For any $\mathbf{I}[[S]]$-algebra $A$, we denote $R_{A,Eis(\mathcal{F},\eta)}$ the local component of $\mathbf{h}(N) \otimes_\iota \mathbf{I}[[S]] \otimes A$ corresponding to the maximal ideal containing $I_{Eis(\mathcal{F},\eta)} \otimes A$ and by $R'_{A,Eis(\mathcal{F},\eta)}$ the $A$-module $R_{A,Eis(\mathcal{F},\eta)}$ divided by its $A$-torsion. Then we define $Eis_A(\mathcal{F},\eta)$ as the kernel of the canonical surjective homomorphism:

$$A \longrightarrow \frac{R'_{Eis(A,\mathcal{F},\eta)}}{R'_{A,Eis(\mathcal{F},\eta)} \cap I_{Eis(\mathcal{F},\eta)} \otimes A}.$$

$Eis_A(\mathcal{F},\eta)$ is what we call the Eisenstein ideal associated to $A$, $\mathcal{F}$ and $\eta$.

## 3 Galois Representations

### 3.1 Galois representation for cuspidal representation of $GSp_4/\mathbf{Q}$

In this section, we recall the main result of R. Weissauer about Galois representation associated with cohomological automorphic representation of $GSp_4(\mathbf{A})$ and we prove under hypothesis over the residual representation some results dealing with the irreducibility of such Galois representation. Let $G_{\mathbf{Q}} = Gal(\bar{\mathbf{Q}}/\mathbf{Q})$.

**Theorem 3.1** *R. Weissauer [35] Let $\pi$ be a cohomological cuspidal representation of weight $(m,n)$ and level $N$. There exists a continuous Galois representation:*
$$\rho_\pi : G_{\mathbf{Q}} \longrightarrow GL_4(\bar{\mathbf{Q}}_p)$$
*unramified at prime not dividing $Np$ such that for all $\ell \nmid Np$, the characteristic polynomial of $\rho_\pi(Frob_\ell)$ is given by $\lambda_\pi(\mathcal{Q}_\ell(X))$. Moreover if $m(\pi) = 1$, then this representation respects a skew-symmetric bilinear form (i.e. $\rho_\pi$ takes values in $GSp_4(\bar{\mathbf{Q}}_p)$.*

**Local properties of $\rho_\pi$.** Let $\rho_{\pi,p}$ be the restriction of $\rho_\pi$ to any decomposition group $D_p$ at $p$. The only known fact about local properties of $\rho_\pi$ at $p$ is the following proposition resulting from works of Faltings, Chai-Faltings and the construction of Weissauer.

**Proposition 3.1** *If $\pi$ is unramified at $p$, then $\rho_{\pi,p}$ is crystalline.*

(i) *If $\pi$ has cohomological weight $(m,n)$, the Hodge-Tate weight of $\rho_{\pi,p}$ are $\{m+n+3, m+2, n+1, 0\}$.*

(ii) *The slopes of the Crystalline Frobenius acting on the filtered $\phi$-module of the local representation $\rho_{\pi,p}$ are contained in the set of p-adic valuations of the roots of $\lambda_\pi(\mathcal{Q}_p(X))$.*

**Proof:** We just give indications of the proof and will give it in details in a subsequent paper. By Weissauer's construction, $\rho_\pi$ is realized in the etale cohomology of the Siegel variety $X_N$ of a level prime $N$ to $p$. By works of Chai-Faltings [2], $X_N$ has smooth compactifications defined over $\mathbf{Z}[\zeta_N, 1/N]$ whose boundary is a normal crossing divisor. Therefore by the main result of Faltings in [5], the representation of $D_p$ on the $p$-adic etale cohomology with coefficient in $L_{m,n}(\mathbf{Q}_p)$ is crystalline with Hodge-Tate weights $(m+n+3, m+2, n+1, 0)$. The result follows from the fact that a subrepresentation of a crystalline representation is crystalline.

The last point comes from the Hecke equivariance of the etale-crystalline comparison isomorphism proven by the Faltings-Jordan's arguments (cf. [6])and the Eichler-Shimura relation on the crystalline side. Note that all the tools needed in the generalization of their proof are given in the Siegel variety case by theorem VI-1.1 of [2]■

**Remark.** When the representation $\pi$ is not a weak endoscopic lift, the fact that the four Hodge-Tate weights occur is equivalent to the stability of the $L$-packet at infinity (i.e. $\pi_f \otimes \pi_\infty^H$ is automorphic if and only if $\pi_f \otimes \pi_\infty^W$ is automorphic too). Moreover, if the Zariski closure of the image of $G_\mathbf{Q}$ is of rank 3, one can see that all Hodge-Tate weights occur thanks to Sen's theory (cf. corollary 1 of [25]).

**Corollary 3.1** *Let $\pi$ be a cohomological cuspidal representation of weight $(m,n)$ which is unramified and ordinary at $p$. We assume that Hodge-Tate weights of $\rho_{\pi,p}$ are $\{m+n+3, m+2, n+1, 0\}$. Then either*
*(i) $\rho_{\pi,p}$ is ordinary at $p$ (i.e.*

$$\rho_\pi|_{I_p} \sim \begin{pmatrix} \chi_p^{m+n+3} & \times & \times & \times \\ 0 & \chi_p^{m+2} & \times & \times \\ 0 & 0 & \chi_p^{n+1} & \times \\ 0 & 0 & 0 & 1 \end{pmatrix}$$

*where $\chi$ denotes the cyclotomic character.)*
*or*
*(ii) the slopes of filtered module of $\rho_{\pi,p}$ are exactly $\{m+2, n+1\}$. In particular, if $m+n+3 < p-1$ then $\rho_{\pi,p}$ is ordinary.*

**Proof:** We know by proposition 3.2 of [28], that $p$-adic valuations of the roots of $\lambda_\pi(\mathcal{Q}_p(X))$ are exactly $\{m+n+3, m+2, n+1, 0\}$. So by proposition above the slopes are contained in $\{m+n+3, m+2, n+1, 0\}$. By autoduality of the representation, one sees that if $\alpha$ is a slope then $a+b+3-\alpha$ is a slope too with same multiplicity. Therefore we have the following possibilities for the slopes:

- $\{m+n+3, m+2, n+1, 0\}$ and the representation is ordinary.

- $\{m+n+3, 0\}$. But this case is impossible because it would imply that the Newton polygon is under the Hodge polygon.

- $\{m+2, n+1\}$ which is the last case we do not know yet how to exclude.

Anyway, the slopes are integral and this implies by the Fontaine-Lafaille'stheory that $\rho_{\pi,p}$ is ordinary when $m+n+3 < p-1$.

**Conjecture 3.1** *If $\pi$ is unramified and ordinary at $p$, then $\rho_{\pi,p}$ is ordinary (i.e. satisfies (3.1.(i))).*

For any irreducible component **J** of the universal p-ordinary Hecke algebra, by corollary 2.2, theorem 3.1 and example 1 of [24] there exist a finite extension $F'$ of the field of fraction $F_\mathbf{J}$ of $\mathbf{J}$ and a unique semi-simple Galois Representation $\rho_\mathbf{J}$ in $GL_4(F')$ unramified outside $Np$ and such that the characteristic polynomial of $Frob_\ell$ is given by :

$$\lambda_\mathbf{J}(\mathcal{Q}_\ell(X))$$

where $\lambda_\mathbf{J}$ is the character of the Hecke algebra corresponding to $\mathbf{J}$. By the results of [28], we get the local property at $p$ of that representation:

**Theorem 3.2** *Let $\mathbf{J}$ of nebentypus $\omega_{a,b}$. Assume that conjecture 3.1 is satisfied, then*

$$\rho_\mathbf{J}|_{I_p} \sim \begin{pmatrix} \xi_1 & \times & \times & \times \\ 0 & \xi_2 & \times & \times \\ 0 & 0 & \xi_3 & \times \\ 0 & 0 & 0 & 1 \end{pmatrix}$$

*with*

$$\xi_1 = (\omega \circ \chi_p)^{a+b+3} <\chi_p>_{u^3(1+T_1)(1+T_2)-1}$$
$$\xi_2 = (\omega \circ \chi_p)^{a+2} <\chi_p>_{u^2(1+T_1)-1}$$
$$\xi_3 = (\omega \circ \chi_p)^{b+1} <\chi_p>_{u(1+T_2)-1}$$

## 3.2 Irreducibility of Galois representations

**Definition 3.1** *A cuspidal representation $\pi$ of $GSp_4$ is called a CAP (Cuspidal Associated to Parabolic) representation if there exist $P$ a proper parabolic subgroup of $GSp_4$ and $\sigma$ a cuspidal representation of the Levi subgroup of $P$ such that $\pi$ and $Ind_{P(\mathbf{A})}^{GSp_4(\mathbf{A})} \sigma$ have the same Langlands parameters at almost all primes.*

The CAP representations of $GSp_4$ are well understood thanks to the works of Piatetski-Shapiro and Soudry. Especially, these authors proved (explicitly) that the CAP representations are in the image of theta correspondence. From their works, one can deduce easily by calculation of the theta lifting at infinity the following proposition. It is also a direct consequence of Theorem 2.5.6 of Harris (cf. [10]).

**Proposition 3.2** *Let $\pi$ be a cohomological representation of weight $(m, n)$. Assume that $\pi$ is a CAP associated to the Klingen parabolic subgroup (resp. to the Siegel parabolic subgroup). Then $n = 0$ (resp. $m = n$).*

**Definition 3.2** *A cuspidal representation $\pi$ of $GSp_4$ is called a weak endoscopic lift, if there exist two cuspidal representations $\sigma_1$, $\sigma_2$ of $GL_2$ with same central character and such that the Langlands'parameters of $\pi$ are the union of those of $\sigma_1$ and $\sigma_2$ at almost all primes.*

It is clear that the representation $\rho_\pi$ is reducible in the CAP and endoscopic case. It is conjectured that outside these two cases the representation $\rho_\pi$ is absolutely irreducible. We prove a partial result towards this conjecture in the following. In order to do this we start by stating a very simple (but useful) variant of the discussion of the paragraph 2 in [25].

**Lemma 3.1** *Let $k$ be a field and $G$ be a group. Let $\rho$ be a reducible representation of $G$ in $GL_4(k)$ and let $\nu_\rho$ be a character of $G$ such that the eigenvalues of $\rho(g)$ come in pairs of the form $\{\alpha, \nu_\rho(g)\alpha^{-1}\}$. Let $H_\rho$ be the Zariski envelope of the image of the semisimplification of $\rho$ in $GL_4(k)$, then we are in one of the following cases:*

1. $H_\rho$ embeds in $(GL_2 \times GL_2)^0 = \{(g_1, g_2) \in GL_2 \times GL_2 \;\; det(g_1) = det(g_2)\}$ the endoscopic group of $GSp_4$. If we write $\rho_\rho = \tau_{1,\rho} \oplus \tau_{2,\rho}$ we have $det(\tau_{1,\rho}) = \nu_\rho$ for $i = 1, 2$ and $\tau_{1,\rho}$ is not isomorphic to $\tau_{2,\rho}$.

2. $H_\rho$ embeds in $GL_2 \times \mathbf{G}_m$, we denote by $\tau_\rho$ the composite of the component in $GL_2$ and the standard representation of $GL_2$ and $\xi_\rho$ the non trivial character corresponding to the component of $\rho$ in $\mathbf{G}_m$. Then we have two possibilities:

    (i) $\rho^{ss} = \tau_\rho \oplus \tau_\rho \otimes \nu_\rho det(\tau_\rho)^{-1}$ (the Siegel Parabolic case); here $\xi_\rho = \nu_\rho det(\tau_\rho)^{-1}$.

    (ii) $\rho^{ss} = \tau_\rho \oplus \xi_\rho \oplus \xi_\rho^{-1} \nu_\rho$ with $det(\tau_\rho) = \nu_\rho$ (the Klingen Parabolic case).

3. $H_\rho$ embeds in $GL_2$, we denote by $\tau_\rho$ the corresponding standard representation. Then we have two possibilities

    (i) $\rho^{ss} = \tau_\rho \oplus \tau_\rho$.

    (ii) $\rho^{ss} = \tau_\rho \oplus 1 \oplus det(\tau_\rho)$.

4. $H_\rho$ embeds to $\mathbf{G}_m^r$ for $r = 1, 2, 3$ (the Borel case).

5. $H_\rho$ embeds in $GL_2$, we denote by $\tau_\rho$ the corresponding standard representation then $\rho^{ss} = Sym^2(\tau_\rho) \oplus 1$

Of course, the representation $\rho_\pi$ verifies the conditions of the above lemma if we take for $\nu_{\rho_\pi}$ the Galois character associated to the central character of $\pi$ by Class Field Theory. Moreover, by some arguments taking into account the Hodge-Tate structure of $\rho_\pi$, the case 5 can be excluded (cf. paragraph 2 of [25]).

Let $k = \bar{\mathbf{F}}_p$ with an odd prime $p$. For any representation $\rho$ of $G_\mathbf{Q}$ in $GL_4(k)$ as in Lemma 3.1 such that $H_\rho$ is in cases 1,2 or 3 above, we consider the following modularity condition:

(**Mod** $\bar{\rho}$) *The Galois representation(s) in $GL_2(\bar{\mathbf{F}}_p)$ interfering in $\bar{\rho}$ are modular (in the sense of Serre) and irreducible on $Gal(\bar{\mathbf{Q}}/\mathbf{Q}(\sqrt{p^*}))$.*

By continuity of $\rho_\pi$, there exists a Galois stable lattice $L$ in the 4-dimensional $\bar{\mathbf{Q}}_p$-space of $\rho_\pi$. Therefore we can consider the residual representation
$$\bar{\rho}_\pi : G_\mathbf{Q} \longrightarrow GL_4(\bar{\mathbf{F}}_p)$$
given by the Galois action on $L \otimes \bar{\mathbf{F}}_p$. If $\bar{\rho}_\pi$ was known to be irreducible, $\rho_\pi$ would be trivially irreducible. Hence we are interesting in the case where $\bar{\rho}$ is reducible.

**Theorem 3.3** *Let $\pi$ be a cuspidal cohomological representation of $GSp_4/\mathbf{Q}$ with central character $\omega_\pi$ such that $\rho_\pi$ is reducible. Assume (**Mod** $\bar{\rho}_\pi$). If $\rho_{\pi,p}$ is ordinary and if $\bar{\rho}_\pi$ is not in case 4, then $\pi$ is a weak endoscopic lift or a CAP representation.*

*In fact we can be a little more precise. If $\bar{\rho}_\pi$ is in case 1 (resp. in case 2), then $\pi$ is a weak endoscopic lift (resp. is a CAP representation associated to the Siegel parabolic subgroup in case 2.(i) and to the Klingen parabolic subgroup in the case 2.(ii)).*

**Proof:** Firstly let us assume that $\rho_\pi$ is in the case 2.(i) or 3.(i); we thus have $\tau_{\bar{\rho}} \cong \bar{\tau}_\rho$. Then obviously $\bar{\rho}_\pi$ will be in the case 2.(i) or (iii). By our assumption, $\bar{\tau}_\rho$ is modular and $\tau_\rho$ is ordinary (as a subrepresentation of an ordinary representation. By Taylor-Wiles'theorem (cf. [37], [26]) and its improvement by F.Diamond (cf. Theorem 5.3 of [4]), we deduce that $\tau_\rho$ is modular; hence is associated with a cuspidal representation $\sigma$ of $GL_{2/\mathbf{Q}}$. Therefore, we see from proposition 2.1 that $\pi$ and $Ind_{Q(\mathbf{A})}^{GSp_4(\mathbf{A})} \sigma \otimes \omega_\pi |.|_\mathbf{A}^2$ have the same Langlands'parameters at almost all primes i.e. $\pi$ is a CAP representation for the Klingen parabolic subgroup. For $\rho_\pi$ in the case 2.(ii) or 3.(ii) cases, we would get by the same way that $\pi$ is CAP for the Siegel parabolic subgroup.

Let us assume now that $\rho_\pi$ is in case 1. By the same sort of arguments, we see then that $\tau_{1,\rho_\pi}$ and $\tau_{2,\rho_\pi}$ are modular and thus $\pi$ is a weak endoscopic lift.

Now we are viewing the situation for $\bar{\rho}$: If $\bar{\rho}_\pi$ is in the case 2.(i), then $\rho_\pi$ is in case 2.(i) and $\pi$ is CAP by the previous discussion. If $\bar{\rho}_\pi$ is in the case 3.(i), then $\rho_\pi$ is in case 2.(i), 3.(i) or 1. In the two first case, we deduce again that $\pi$ is CAP, but in the third (case 1), $\pi$ is a weak endoscopic lift. The other cases, can be deduced similarly. ∎

**Definition 3.3** *An irreducible closed subscheme $Spec(\mathbf{J})$ is called endoscopic if it contains a densely populated set of arithmetic points associated to endoscopic representations.*

**Remarks:** Assume that **J** is an irreducible component of $\mathbf{h}(N)$. If $N = 1$, by using the Taylor-Wiles'theorem (with same hypothesis as theorem 3.3), one can prove that there exist $\mathbf{I}_1$ and $\mathbf{I}_2$ some finite extensions of $O_K[[T]]$ and $\mathbf{I}_i$-adic ordinary elliptic cusp form $F_i$ for i=1 and 2, such that:

$$\lambda_{\mathbf{J}}(\mathcal{Q}_\ell(X)) = (X^2 - a(\ell; F_1)X - \ell\lambda_{F_1}(T(\ell,\ell)))((X^2 - a(\ell; F_2)X - \ell\lambda_{F_2}(T(\ell,\ell))))$$

In general, we are confronted to the problem of bounding the level of the couple of elliptic modular forms associated to an endoscopic lift. However, if $\rho_{\mathbf{J}}$ is reducible one can easily prove that $\rho_{\mathbf{J}} = \rho_1 \oplus \rho_2$ with $det(\rho_1) = det(\rho_2)$.

**Theorem 3.4** *Assume conjecture 3.1. Let $Spec(\mathbf{J})$ be a Klingen-regular (resp. Siegel-regular) irreducible closed subscheme of $Spec(\mathbf{h}(N))$. If either $\bar\rho_{\mathbf{J}}$ is in case 2.(i) (resp. 2.(ii)) or $\bar\rho_{\mathbf{J}}$ is in case 3.(i) (resp. 3.(ii)) and $\mathbf{J}$ is not endoscopic.*

*Then $\rho_{\mathbf{J}}$ is absolutely irreducible.*

**Proof** Let us fix a stable lattice $L \in F'^4$. Let $\mathbf{J}'$ the integral closure of $\mathbf{J}$ in $F'$. Since $\mathbf{J}$ is not endoscopic, there exists an arithmetic dominant and ?-regular weight $\psi_0$ such that $\lambda_{\mathbf{J}}$ mod $P_{\psi_0}$ is not endoscopic and $P'$ a prime of $\mathbf{J}'$ over $P_{\psi_0}$ such that $L_{P'}$ is free of rank 4. Then we can consider $\rho_{\mathbf{J}}$ mod $P'$ in $GL_4(\mathbf{J}'_{P'}/P')$ that is isomorphic to $\rho_{\lambda_{\mathbf{J}} \, mod \, P_{\psi_0}}$. On the other hand, by Proposition 3.2, $\lambda_{\mathbf{J}}$ mod $P_{\psi_0}$ is not CAP because $\psi_0$ is regular. By theorem 3.2, $\rho_{\lambda_{\mathbf{J}} mod \, P'}$ is ordinary, therefore it is absolutely irreducible by theorem 3.3. Therefore $\rho_{\mathbf{J}}$ is too. ∎

We now discuss about the following conjecture:

**Conjecture 3.2** *For $\pi$ cuspidal and cohomological, $\rho_\pi$ respects (up to a multiplicative factor) a skew-symmetric form.*

If we assume multiplicity one holds for a cuspidal representation $\pi$ of $GSp_4$, this conjecture follows from the realization of these representations in the degre 3 etale cohomology of the Siegel threefold, the isomorphism $\pi \cong \pi^* \otimes \nu_\pi^{-1}$ and the existence of the cup-product. Therefore, if a representation is quasi-equivalent (at all but a finite set of place) to a generic representation then $\rho_\pi$ preserves a symplectic form. In general, we can just state the following:

**Proposition 3.3** *$\rho_{\mathbf{J}}$ is self dual (up to a twist). Moreover, if one irreducible specialization of $\rho_{\mathbf{J}}$ takes values in $GSp_4$ or if $\rho_{\mathbf{J}}$ is reducible, then $\rho_{\mathbf{J}}$ takes values in $GSp_4$.*

**Proof:** Consider $\nu_{\mathbf{J}}$ the Galois character in $\mathbf{J}^\times$ unramified outside $Np$ such that $\nu_{\mathbf{J}}(Frob_\ell) = \ell^3 \lambda_{\mathbf{J}}(S_\ell)$. Then we can see that by theorem 3.1 and corollary 2.2, we have $tr(\rho_{\mathbf{J}}(\sigma)) \equiv tr(\check{\rho}_{\mathbf{J}} \otimes \nu_{\mathbf{J}}(\sigma))$ mod $P_\psi$ for all $\sigma \in G_{\mathbf{Q}}$ and dominant regular $\psi$. Therefore by Brauer-Nesbitt's theorem, $\rho_{\mathbf{J}} \cong \check{\rho}_{\mathbf{J}} \otimes \nu_{\mathbf{J}}$ and therefore $\rho_{\mathbf{J}}$ takes values in a similitude group for some bilinear form $B$. If $\rho_{\mathbf{J}}$ is irreducible, a standard argument proves that $B$ is either symmetric or skew-symmetric. If $B$ mod $P_\psi$ is skew symmetric for one $P_\psi$, $B$ is thus skew-symmetric too. The reducible case is obvious and follows from the classification of lemma 3.1.■

## 3.3 The Selmer groups

We recall now the definition of the Selmer groups we are interested in. We use the definitions and notations of section 2.4. Let $\rho_{\mathcal{F}}$ be the Galois representation in $GL_2(\mathbf{I})$ associated to $\mathcal{F}$; it is continuous and unramified outside $Np$. For each $\ell \nmid$ the characteristic polynomial of $Frob_\ell$ is given by

$$X^2 - a(\ell; \mathcal{F})X + \ell^{-1}\epsilon_{\mathcal{F}}(\ell) <\ell>_T$$

Moreover $\rho_{\mathcal{F}}$ is ordinary at $p$. Then if we assume the hypothesis:

(**Reg** F) $(\epsilon_{\mathcal{F}})_p \neq \omega$

there exists $g_+ \in GL_2(\mathbf{I})$ such that for all $\sigma \in I_p$

$$\rho_{\mathcal{F}}(\sigma) = g_+ \begin{pmatrix} det(\rho_{\mathcal{F}})(\sigma) & \star \\ 0 & 1 \end{pmatrix} g_+^{-1}$$

Let us consider now the Galois character in $O_K[[S]]$ defined by:

$$\tilde{\eta}(\sigma) = \eta^G(\sigma) <\chi_p(\sigma)>_S$$

where $\eta^G$ denotes the Galois character associated to the Dirichlet character $\eta$ by Class Field Theory; obviously, $\tilde{\eta}$ is the universal deformation of the character $\eta^G$. For any Dirichlet character $\psi$, we will denote by $\psi_p$ its component on the factor of conductor a power of $p$.

Let $ad^0(\rho_{\mathcal{F}}) \otimes \tilde{\eta}$ the adjoint representation of $G_{\mathbf{Q}}$ on $sl_2(\mathbf{I}) \otimes \mathbf{Z}_p[[S]]$ twisted by $\tilde{\eta}$. We set $F^+(ad^0(\rho_{\mathcal{F}}))$ the sub $\mathbf{I}$-module of $sl_2(\mathbf{I})$ generated by the matrices of the following type:

$$g_+ \begin{pmatrix} 0 & \star \\ 0 & 0 \end{pmatrix} g_+^{-1}$$

By ordinarity of $\rho_{\mathcal{F}}$, this submodule is stable by the action of the inertia $I_p$, therefore we can consider the action of $I_p$ on $ad^0(\rho_{\mathcal{F}})/F^+(ad^0(\rho_{\mathcal{F}})) \otimes \tilde{\eta}$. We fix $I_q$ an inertia subgroup for all prime $q$. For any irreducible quotient $A$ of $\mathbf{I}[[S]]$ and any finite set of primes $\Sigma$ containing $p$, we define $Sel_{A,\Sigma}(ad^0(\rho_{\mathcal{F}}) \otimes \tilde{\eta})$ as the kernel of

$$H^1(G_{\mathbf{Q}}, (ad^0(\rho_{\mathcal{F}}) \otimes \tilde{\eta}^{-1}) \otimes A^*) \longrightarrow \bigoplus_{q \notin \Sigma} H^1(I_q; (ad^0(\rho_{\mathcal{F}}) \otimes \tilde{\eta}^{-1}) \otimes A^*)$$
$$\oplus H^1(I_p; ((ad^0(\rho_{\mathcal{F}})/F^+(ad^0(\rho_{\mathcal{F}}) \otimes \tilde{\eta}^{-1}) \otimes A^*)$$

with $A^* = Hom_{\mathbf{Z}_p}(A, \mathbf{Q}_p/\mathbf{Z}_p)$ the Pontrjagin dual of $A$.

We consider the following condition on A:

(**Arith** A) *A is an arithmetic Klingen-regular irreducible quotient of* $\mathbf{I}[[S]]$ *(cf. definition 2.2)*

We denote by $\theta_A$ the corresponding surjective homomorphism $\mathbf{I}[[S]] \to A$. Note that $A$ is still a Krull ring.

This condition is fulfilled in the following interesting cases:

1. $A = \mathbf{I}[[S]]$ the Selmer group is the two variable one considered in [19].

2. $A = O_K[[S]]$ where the map $\theta_A$ is obtained from a modular form $f_k$ of weight $k > 3$ belonging to the family $\mathcal{F}$. In that case, the Selmer group is the classical one-variable Selmer group attached to the cyclotomic twists of the symmetric square of the modular representation $\rho_{f_k}$.

3. $A = \mathbf{I}$ with $\theta_A(S) = (1+u)^l - 1$ with $l > 1$. Then the Selmer group has only the weight-variable of the family $\mathcal{F}$.

4. $A = O_K$ where $\theta_A(S) = (1+u)^l - 1$ with $k - 3 \geq l > 1$ and $\theta_A$ is obtained from a modular form $f_k$ of weight $k > 3$ belonging to the family $\mathcal{F}$.

**Remark:** The case $l = 0$ and $k = 2$ could be also treated by the same way. We have not included it here because it should be studied by using directly the Hecke algebra of $GSp_4$ associated to automorphic forms whose archimedean component is the holomorphic completely degenerated limit of dicrete series (because it is not controled by our universal Hecke algebra). Note that the existence of the Galois representations for such forms can be constructed by the way of congruences with higher weights forms by using the Hasse invariants constructed by Blasius-Ramakrishnan and Clozel.

A special case of one conjecture of [8] is

**Conjecture 3.3** $Sel_{A,\Sigma}(ad^0(\rho_{\mathcal{F}}) \otimes \tilde{\eta}^{-1})$ *is co-torsion over* $A$. *We denote by* $F^{\Sigma}_{ad^0(\rho_{\mathcal{F}}) \otimes \tilde{\eta}^{-1}, A}$ *its characteristic ideal in* $A$.

**Remark:** This conjecture is proven by Hida for $A = \mathbf{I}[[S]]$ when $\eta$ is trivial cf. [15] or for even $\eta$ ( which are unramified at $p$) when $\Sigma$ is reduced to

{p} see [16] and [17]. It seems that his method should solve the general case provided that the isomorphisms betweenthe Hecke algebra over the totally real field fixed by the kernel of $\eta$ and the corresponding universal deformation ring is proved in the unrestricted cases (see [7], [16] and [17]).

Let $\Sigma_N$ be the set of primes dividing $N$ and $A$ satisfying (**Reg** $A$). Our main result is:

**Theorem 3.5** *Let $\mathcal{F}$ be an **I**-adic cusp form with $\epsilon_{\mathcal{F}} \neq \omega$. We assume conjecture 3.1 and either conjecture 3.2 or that $L_p(\eta^{-1})$ the Kuboto-leopoldt p-adic L function associated to the branch $\eta^{-1}$ is a unit in $O_K[[S]]$. Let $P$ be an height one prime ideal of $A$ such that*

(i) $P \neq (\theta_A(S))$ *if $\eta_p$ is trivial modulo $p$.*

(ii) $P \cap \theta_A(\mathbf{Z}_p[[T,S]]) \neq (\theta_A((1+S)(1+T)-u))$ *if $\eta_p(\epsilon_{\mathcal{F}})_p$ is trivial modulo $p$.*

*then*
$$\mathrm{length}_{A_P}(Sel_{A,\Sigma_N}(ad^0(\rho_{\mathcal{F}}) \otimes \tilde{\eta}^{-1})_P) \geq v_P(Eis_A(\mathcal{F},\eta)).$$

*Therefore, if we assume conjecture 3.3 for $Sel_{A,\Sigma_N}(ad^0(\rho_{\mathcal{F}}) \otimes \tilde{\eta}^{-1})$, we have the divisibility*
$$Eis_A(\mathcal{F},\eta) | F^{\Sigma_N}_{ad^0(\rho_{\mathcal{F}}) \otimes \tilde{\eta}^{-1},A}$$

*in $A$ (resp. in $A_{\theta_A(S)}$ if $\eta_p$ is trivial modulo $p$ and in $A_{\theta_A((1+S)(1+T)-u)}$ if $\eta_p(\epsilon_{\mathcal{F}})_p$ is trivial modulo $p$).*

**Remarks:** a) Note that $Eis_A(\mathcal{F},\eta)$ is not a divisor and the divisibility $Eis_A(\mathcal{F},\eta) | F^{\Sigma_N}_{ad^0(\rho_{\mathcal{F}}) \otimes \tilde{\eta}^{-1},A}$ means that the divisor associated to $Eis_A(\mathcal{F},\eta)$ divides $F^{\Sigma_N}_{ad^0(\rho_{\mathcal{F}}) \otimes \tilde{\eta}^{-1},A}$.

b) In case (i) for $A = \mathbf{I}[[S]]$, by results of Greenberg-Tilouine (cf. [9]) and Hida (cf [15]), the $(S)$-adic valuations of the $p$-adic L function $L_{ad^0(\rho_{\mathcal{F}}) \otimes \tilde{\eta}^{-1}}$ and the characteristic ideal $F^{\Sigma_N}_{ad^0(\rho_{\mathcal{F}}) \otimes \tilde{\eta}^{-1}, \mathbf{I}[[S]]}$ are both equal to 1 (see also [19]). In case (ii), it is easy to verify that $v_P(L_{ad^0(\rho_{\mathcal{F}}) \otimes \tilde{\eta}^{-1}}) = 0$ using the interpolation property of this $p$-adic L function. Therefore these assumptions are not restrictions in the proof of the divisibility $L_{ad^0(\rho_{\mathcal{F}}) \otimes \tilde{\eta}^{-1}} \mid F^{\Sigma_N}_{ad^0(\rho_{\mathcal{F}}) \otimes \tilde{\eta}^{-1}}$.

### 3.4 Proof of theorem 3.5

We fix once for all, a finite extension $L$ of the field of fractions $F_A$ of $A$ such that the Galois representations associated to the irreducible components of $R'_{A,\mathcal{F},\eta}$ take values in $GSp_4(L)$. Since $R'_{A,\mathcal{F},\eta} \otimes L$ is semi-simple by corollary 2.4, it is possible to glue these Galois representations and thus to get a Galois action on $V_L = (R'_{A,\mathcal{F},\eta} \otimes L)^4$. Let $A_L$ be the integral closure of $A$ in $L$. From now one we fix $P$ a height one prime of $A$ dividing $Eis_A(\mathcal{F},\eta)$ as in the theorem and $Q$ a height one prime of $A_L$ lying over $P$. We denote by $A_Q$ the localization of $A_L$ at $Q$. It is a discrete valuation ring for which we choose a uniformizing parameter $\varpi_Q$ and we denote by $e$ the ramification index of $P$ in $A_L$ (i.e $P.A_Q = \varpi_Q^e A_Q$).

Since $P$ divides $Eis_A(\mathcal{F},\eta)$, by our assumption (Reg A), one sees from Theorem 3.4 that the Galois representation on $V_L$ is irreducible as $R_{A,\mathcal{F},\eta} \otimes L$-module. Indeed, let us take $\mathbf{J}$ be any irreducible component of $R'_{A,\mathcal{F},\eta}$ occurring in $R'_{A,\mathcal{F},\eta} \otimes A_Q$. Then it follows from the congruence relations that the characteristic polynomial of $\rho_{\lambda_\mathbf{J}} \otimes 1_A$ is congruent to that of $(\rho_\mathcal{F} \oplus \rho_\mathcal{F} \otimes \tilde{\eta}) \otimes 1_A$ modulo $P$. Now if $\rho_{\lambda_\mathbf{J}}$ was endoscopic and reducible, its semi-simplification would be isomorphic to $\rho_1 \oplus \rho_2$ with $det(\rho_1) = det(\rho_2)$ but this implies that $\theta_A(\tilde{\eta})$ is trivial modulo $P$ what is impossible by assumption (i) of the theorem. Therefore, $\rho_{\lambda_\mathbf{J}}$ is irreducible; moreover its image falls in $GSp_4(L)$ by proposition 3.3 and corollary **??**. We can thus apply section 1.3 to our situation by taking $\mathcal{B} = R'_{A,\mathcal{F},\eta} \otimes A_Q$, $I = I_{Eis(\mathcal{F},\eta)} \otimes A_Q$, $\rho_0 = \rho_\mathcal{F} \otimes 1_{A_Q}$, $\nu_0 = \tilde{\eta} det(\rho_\mathcal{F}) \otimes 1_{A_Q}$. If we denote by $n = v_P(Eis_A(\mathcal{F},\eta))$, then we have
$$\mathcal{B}/I = R'_{A,\mathcal{F},\eta} \otimes A_Q / I_{Eis(\mathcal{F},\eta)} \otimes A_Q = A_Q/Q^{ne}$$
and thus there exists a $A_Q$-lattice in $V_L$ without quotient isomorphic to $\rho_\mathcal{F} \otimes 1_{A_Q}$ mod $P$ and such that its reduction modulo $I$ provides us the following exact sequence:

$$0 \longrightarrow V(\rho_\mathcal{F}) \otimes N \longrightarrow \mathcal{L} \otimes_{\lambda_{\mathcal{F},\eta}} A_Q/Q^{ne} \xrightarrow{s} V({}^t\rho_\mathcal{F}^{-1} \otimes det(\rho_\mathcal{F})\tilde{\eta}) \otimes A_Q/Q^{ne} \longrightarrow 0$$

where

1. $N = \mathcal{T}/I\mathcal{T}$ with $\mathcal{T}$ faithful $A_Q$-module (contained in the total ring of fraction of $R'_{A,\mathcal{F},\eta} \otimes A_Q$).

2. $s$ is a section of $A_Q$-module but is not Galois equivariant.

3. $\mathcal{L} \otimes_{\lambda_{\mathcal{F},\eta}} A_Q/Q^{ne}$ is endowed with a skew symmetric bilinear form such that
$$\bar{\phi}(\begin{pmatrix} t_1 \\ t_2 \end{pmatrix}, \begin{pmatrix} s(b_1) \\ s(b_2) \end{pmatrix}) = b_1.t_1 + b_2.t_2$$
for all $\begin{pmatrix} t_1 \\ t_2 \end{pmatrix} \in V(\rho_\mathcal{F}) \otimes N$ and $\begin{pmatrix} b_1 \\ b_2 \end{pmatrix} \in V({}^t\rho_\mathcal{F}^{-1}) \otimes A_Q/Q^{ne}$.

Note that by point 1, one sees easily by proposition 4 and 7 of the appendix of [22] that the
$$length_{A_Q}(N) \geq ne$$

Consider now the cocycle $C_{Eis}$ from $G_{\mathbf{Q}}$ in $M_4(A_Q) \otimes_{A_Q} N$ defined by

$$C_{Eis}(g)(m) = m - g.s(g^{-1}.m)$$

for all $m \in V({}^t\rho_{\mathcal{F}}^{-1} \otimes det(\rho_{\mathcal{F}})\tilde{\eta}) \otimes A_Q/Q^{ne}$ and where we have done the identification

$$Hom_{A_Q}(V({}^t\rho_{\mathcal{F}}^{-1} \otimes det(\rho_{\mathcal{F}})\tilde{\eta}) \otimes A_Q/Q^{ne}, V(\rho_{\mathcal{F}})) \otimes N \cong M_2(A_Q) \otimes_{A_Q} N$$

by using the canonical basis of $V(\rho_{\mathcal{F}})$ and $V({}^t\rho_{\mathcal{F}}^{-1}) \otimes det(\rho_{\mathcal{F}})\tilde{\eta}$ Let us note that the action of $\gamma \in G_{\mathbf{Q}}$ on $X \otimes n \in M_2(A_Q) \otimes_{A_Q} N$ induced by this identification is given by:

$$g.(X \otimes n) = (det(\rho_{\mathcal{F}})\tilde{\eta}(g))^{-1}.\rho_{\mathcal{F}}(g) X {}^t\rho_{\mathcal{F}}(g) \otimes n$$

Let us consider the map:

$$Hom_{A_Q}(N; Q^{-ne}A_Q/A_Q) \xrightarrow{\iota_N} H^1(G_{\mathbf{Q}}; M_2(A_Q) \otimes_{A_Q} Q^{-ne} A_Q/A_Q)$$

where $\iota_N(f)$ is defined as the cohomology class of $C_f = (1 \otimes f) \circ C_{Eis}$.

**Fact 1.** $\iota_N$ *is injective.*

**Proof:** Indeed assume a moment it is not and let $f \in Ker(\iota_N)$. If $f \neq 0$, there exist $N'$ with $N \supset N' \supset Ker(f)$ and $N/N' \cong A_Q/Q$. And let us set

$$M = \frac{\mathcal{L} \otimes_{\lambda_{\mathcal{F},\eta}} A_Q/Q^{ne}}{Q.\mathcal{L} + V(\rho_{\mathcal{F}}) \otimes N'}$$

Then since $\iota_N(f) = 0$, we see easily that $M$ as $G_{\mathbf{Q}}$-module is isomorphic to

$$\left(V(\rho_{\mathcal{F}}) \oplus V({}^t\rho_{\mathcal{F}}^{-1} \otimes det(\rho_{\mathcal{F}})\tilde{\eta})\right) \otimes A_Q/Q.$$

But this contradicts the fact that $\mathcal{L}$ does not have a quotient isomorphic to $V(\rho_{\mathcal{F}}) \otimes A_Q/Q$. ∎

In order to relate the cocycles we have constructed to elements of our Selmer group, let us first remark that the representation of $G_{\mathbf{Q}}$ defined above is isomorphic to $ad(\rho_{\mathcal{F}}) \otimes \tilde{\eta}^{-1}$; we note $\iota'_N$ the composite of $\iota_N$ with this isomorphism.

**Fact 2.** *For all* $f \in Hom_{A_Q}(N; Q^{-ne}A_Q/A_Q)$ *and* $g \in G_{\mathbf{Q}}$, $C'_f(g) = \iota'_N(f)(g) \in ad^0(\rho_{\mathcal{F}}) \otimes \tilde{\eta}^{-1} \otimes Q^{-ne}A_Q/A_Q$.

**Proof:** We have something to prove only when the $L_p(\eta^{-1})$ is not supposed to be a unit in $O_K[[S]]$ because $Sel(ad(\rho_{\mathcal{F}}) \otimes \tilde{\eta}^{-1}) = Sel(ad^0(\rho_{\mathcal{F}}) \otimes \tilde{\eta}^{-1}) \oplus Sel(\tilde{\eta}^{-1})$ and the last factor has characteristic ideal equal to $L_p(\eta^{-1})$ by the Iwasawa conjecture proven by Mazur-Wiles. Let $r > 0$ such that $Im(f) = Q^{-r}A_Q/A_Q$. Then the action of any $g \in G_{\mathbf{Q}}$ on

$$M_f = \frac{\mathcal{L} \otimes_{\lambda_{\mathcal{F},\eta}} A_Q/Q^{ne}}{Q^r.\mathcal{L} + V(\rho_{\mathcal{F}}) \otimes Ker(f)} \cong (A_Q/Q^r A_Q)^4$$

The matrix of the action of $g \in G_{\mathbf{Q}}$ with respect of the basis of $M_f$ constructed from the canonical basis of $V(\rho_{\mathcal{F}}) \otimes N/Ker(f)$ and the image by $s$ of the canonical basis of $V({}^t\rho_{\mathcal{F}}^{-1}) \otimes A_Q/Q^{ne}$ is given by:

$$E(g) = \begin{pmatrix} \rho_{\mathcal{F}}(g) & Q^r.C_f(g){}^t\rho_{\mathcal{F}}(g)^{-1} \\ 0_2 & {}^t\rho_{\mathcal{F}}(g)^{-1}det(\rho_{\mathcal{F}})\tilde{\eta}(g) \end{pmatrix} \text{ mod. } Q^r$$

Since by the point 3 above (cf. proposition 1.1), $E(g) \in GSp_4(A_Q/Q^r A_Q)$ and therefore $C_f(g)$ is symmetric. This implies that $C'_f(g) = \iota'_N(f)(g) \in (ad^0(\rho_{\mathcal{F}}) \otimes \tilde{\eta}^{-1}) \otimes \varpi_Q^{-ne}A_Q/A_Q$. ∎

Since the Galois representations occurring in the above construction are unramified outside $Np$, it follows that it is the same for $C'_f$. Then the local condition is satisfied for $q \nmid Np$.

**Study of the local condition at $p$.**
Let us prove now that $C'_f$ is ordinary in the sense that $C'_f(g)$ modulo $F^+(ad^0(\rho_{\mathcal{F}}) \otimes \tilde{\eta}^{-1}) \otimes \varpi_Q^{-ne}A_Q/A_Q$ is a coboundary. To prove that we need temporary to assume Conjecture 3.1; let us do it. By theorem 3.2, this implies there exists a $(R'_{A,\mathcal{F},\eta} \otimes F_A)$-submodule of rank 1 in $V_L$ which is stable by the inertia group $I_p$ with action given by $\xi_1$ (cf. theorem 3.2). Moreover by assumption on $P$, the characters $(\xi_i)_{1 \leq i \leq 4}$ are different modulo $P$. This implies easily that $\mathcal{L}[\xi_1]$ is a direct factor of $\mathcal{L}$. Therefore there exist $v_0 \in M_f \backslash \varpi_Q M_f$ such that $\sigma.v_o = \tilde{\eta}(\sigma) det(\rho_{\mathcal{F}})(\sigma).v_0$ for all $\sigma \in I_p$. Let $g_+$ such that

$$\rho_{\mathcal{F}}(\sigma) = g_+ \begin{pmatrix} det(\rho_{\mathcal{F}})(\sigma) & t(\sigma) \\ 0 & 1 \end{pmatrix} g_+^{-1}$$

for all $\sigma \in I_p$ and let us write $v_0$ in the basis used above: $v_0 = {}^t(\alpha, \beta, \gamma, \delta)$ and set ${}^t(\gamma', \delta') = {}^tg_+{}^t(\gamma, \delta)$. Then writing matricially the equality $\sigma.v_o = \tilde{\eta}(\sigma)det(\rho_{\mathcal{F}})(\sigma).v_0$, we get $\gamma' = 0$ and

$$\rho_{\mathcal{F}}(\sigma) \begin{pmatrix} \alpha \\ \beta \end{pmatrix} + C(\sigma){}^t\rho_{\mathcal{F}}(\sigma)^{-1}\tilde{\eta}det(\rho_{\mathcal{F}})(\sigma){}^tg_+^{-1}\begin{pmatrix} 0 \\ \delta' \end{pmatrix} = \tilde{\eta}det(\rho_{\mathcal{F}})(\sigma)\begin{pmatrix} \alpha \\ \beta \end{pmatrix}$$

Since $v_0 \notin \varpi_Q M_f$ this implies that $\delta' \not\equiv 0$ mod $Q$ and therefore we can assume $\delta' = 1$. Writing

$$g_+^{-1}C_f(\sigma){}^tg_+^{-1} = \begin{pmatrix} a_\sigma & b_\sigma \\ b_\sigma & c_\sigma \end{pmatrix}$$

we thus have:

$$\begin{pmatrix} b_\sigma \\ c_\sigma \end{pmatrix} = \tilde{\eta} det(\rho_{\mathcal{F}})(\sigma)^{-1} \begin{pmatrix} det(\rho_{\mathcal{F}})(\sigma) & t(\sigma) \\ 0 & 1 \end{pmatrix} g_+^{-1} \begin{pmatrix} \alpha \\ \beta \end{pmatrix} - g_+^{-1} \begin{pmatrix} \alpha \\ \beta \end{pmatrix}$$

A simple computation implies that

$$g_+^{-1} C_f(\sigma)^t g_+^{-1} =$$

$$\begin{pmatrix} \star & 0 \\ 0 & 0 \end{pmatrix} + [det(\rho_{\mathcal{F}}(\sigma))\tilde{\eta}(\sigma)]^{-1} g_+^{-1} \rho_{\mathcal{F}}(\sigma) \begin{pmatrix} 0 & \alpha \\ \alpha & \beta \end{pmatrix} {}^t\rho_{\mathcal{F}}(\sigma)^t g_+^{-1} - g_+^{-1} \begin{pmatrix} 0 & \alpha \\ \alpha & \beta \end{pmatrix} g_+^{-1}$$

Taking into account the Galois equivariant isomorphism

$$ad^0(\rho_{\mathcal{F}}) \otimes \tilde{\eta} \cong Sym^2(\rho_{\mathcal{F}}) \otimes det(\rho_{\mathcal{F}})^{-1} \tilde{\eta}^{-1},$$

we get that $C'_f|_{I_p}$ is trivial in

$$H^1(I_p; (ad^0(\rho_{\mathcal{F}})/F^+(ad^0(\rho_{\mathcal{F}})) \otimes \tilde{\eta}^{-1}) \otimes \varpi_Q^{-ne} A_Q/A_Q).$$

∎

We conclude from the above discussion that the length of $Sel_{A_Q/Q^{en}, \Sigma_N}(ad^0(\rho_{\mathcal{F}}) \otimes \tilde{\eta}^{-1}))$ is greater or equal to $ne$ and therefore

$$length_{A_P}(Sel_{A_P/P^n, \Sigma_N}(ad^0(\rho_{\mathcal{F}}) \otimes \tilde{\eta}^{-1})) \geq n$$

because

$$Sel_{A_P/P^n, \Sigma_N}(ad^0(\rho_{\mathcal{F}}) \otimes \tilde{\eta}^{-1}) \hookrightarrow Sel_{A, \Sigma_N}(ad^0(\rho_{\mathcal{F}}) \otimes \tilde{\eta}^{-1})_P[P^n]$$

by irreducibility of $ad^0(\rho_{\mathcal{F}}) \otimes \tilde{\eta}^{-1} \otimes A_P/P.A_P$.∎

**Remark:** Actually, a similarly proof of the above theorem should give the following:

Let $\mathcal{F}$ be an **I**-adic cusp form with $\epsilon_{\mathcal{F}} \neq \omega$. We assume conjecture 3.1 and either conjecture 3.2 or that $L_p(\eta^{-1})$ the Kuboto-leopoldt p-adic L function associated to the branch $\eta^{-1}$ is a unit in $O_K[[S]]$. If $\eta_p$ and $\eta_p(\epsilon_{\mathcal{F}})_p$ are not trivial modulo p. then

$$\text{Fitt}_A(Sel_{A, \Sigma_N}(ad^0(\rho_{\mathcal{F}}) \otimes \tilde{\eta}^{-1})) \subset Eis_A(\mathcal{F}, \eta).$$

Moreover, if $dim_\kappa((Sel_{\kappa, \Sigma_N}(ad^0(\bar{\rho}_{\mathcal{F}}) \otimes \bar{\eta}^{-1}) = 1$, then we have an injective homomorphism:

$$[A/Eis_A(\mathcal{F}, \eta)]^* \hookrightarrow Sel_{A, \Sigma_N}(ad^0(\rho_{\mathcal{F}}) \otimes \tilde{\eta}^{-1})$$

In light of the above remark, we like to state the following conjecture:

**Conjecture 3.4** *Let $\mathcal{F}$ be an $\mathbf{I}$-adic cusp form with $\epsilon_\mathcal{F} \neq \omega$. If $\eta_p$ and $\eta_p(\epsilon_\mathcal{F})_p$ are not trivial modulo $p$. then*

$$\mathrm{Fitt}_A(Sel_{A,\Sigma_N}(ad^0(\rho_\mathcal{F}) \otimes \tilde{\eta}^{-1})) = Eis_A(\mathcal{F}, \eta).$$

*Especially, if $dim_\kappa((Sel_{\kappa,\Sigma_N}(ad^0(\bar{\rho}_\mathcal{F}) \otimes \bar{\eta}^{-1}) = 1$, then we have an isomorphism:*

$$[A/Eis_A(\mathcal{F}, \eta)]^* \cong Sel_{A,\Sigma_N}(ad^0(\rho_\mathcal{F}) \otimes \tilde{\eta}^{-1})$$

### 3.5 An application to congruences for base change.

We give below a discussion of an application of Theorem 3.5. Let $f$ be an ordinary cuspidal newform of weight $k$ and $\hat{f}$ the base change to $\mathbf{Q}(\eta)$ the totally real field associated to $\eta$ by Class Field Theory. Let us state the following conjecture:

**Conjecture 3.5** *Assume that $\eta$ is even. If $p$ divides*

$$\frac{L(1, Ad(\rho_f) \otimes \eta^{-1})}{\Omega_{Hida}(f)},$$

*then there exists $g$ an Hilbert modular (cuspidal) newform of weight $k$ which is not a base change from $\mathbf{Q}$ and such that*

$$g \equiv \hat{f} \mod \mathfrak{P}$$

*for $\mathfrak{P}$ a prime of $\bar{\mathbf{Q}}$ above $p$.*

Let $K$ be a finite extension of $\mathbf{Q}_p$ containing the eigenvalues of $f$ and the values of $\eta$. Consider $Eis(f, \eta) = \theta_{O_K}(Eis(\mathcal{F}, \eta)$ where $\mathcal{F}$ is an $\mathbf{I}$-adic Hida family whose a specialization $\mathbf{I} \to O_K$ gives $f$. Then we have the following result towards the above conjecture:

**Corollary 3.2** *We assume conjecture 3.1 and either conjecture 3.2 or that $L_p(\eta^{-1})$ the Kuboto-leopoldt $p$-adic L function associated to the branch $\eta^{-1}$ is a unit in $O[[S]]$. Assume that*

- $\mathbf{Q}(\eta)$ *and* $\mathbf{Q}(\mu_p)$ *are linearly disjoint.*

- $\eta_p$ *and* $\eta_p\omega^{k-1}\epsilon_f$ *are non trivial modulo $p$*

- *For all $q|N$, $H^1(I_q, ad^0(\bar{\rho}_f) \otimes \bar{\eta}^{-1}) = 0$.*

- $k > 3$

- $\bar{\rho}_f|_{Gal(\bar{\mathbf{Q}}/\mathbf{Q}(\eta)(\sqrt{p^*}))}$ *is absolutely irreducible*

*Then if $p$ divides $Eis(f, \eta)$, there exist $g$ an Hilbert modular (cuspidal) newform of weight $k$ which is not a base change from $\mathbf{Q}$ and $\mathfrak{P}$ a prime of $\bar{\mathbf{Q}}$ above $p$ such that*

$$g \equiv \hat{f} \mod \mathfrak{P}.$$

**Proof:** By Theorem 3.5 for $A = O_K$ and $l = 2$ in case 4 of the list of examples given before the theorem, we see that

$$Sel_{O_K, \Sigma_N}(ad^0(\rho_f) \otimes \eta^{-1} <\chi_p>^{-2}) \neq 0.$$

Moreover, by the second hypothesis, we see that:

$$Sel_{O_K, \Sigma_N}(ad^0(\rho_f) \otimes \eta^{-1} <\chi_p>^{-2}) = Sel_{O_K, \emptyset}(ad^0(\rho_f) \otimes \eta^{-1} <\chi_p>^{-2})$$

Let us denote by $\mathbf{F}$ the residual field of $O_K$, then we have by irreducibility of $ad^0(\bar{\rho}_f)$:

$$Sel_{O_K, \emptyset}(ad^0(\rho_f) \otimes \eta^{-1} <\chi_p>^{-2})[p] =$$
$$= Sel_{\mathbf{F}, \emptyset}(ad^0(\bar{\rho}_f) \otimes \bar{\eta}^{-1}) = Sel_{O_K, \emptyset}(ad^0(\rho_f) \otimes \eta^{-1})[p].$$

We deduce that

$$(\star) \quad Sel_{O_K, \emptyset}(ad^0(\rho_f) \otimes \eta^{-1}) \neq 0.$$

We can also consider the Selmer group $Sel_{O_K, \emptyset}^{\mathbf{Q}(\eta)}(ad^0(\rho_f))$ over the totally real field $\mathbf{Q}(\eta)$ defined by considering the action restricted to $Gal(\bar{\mathbf{Q}}/\mathbf{Q}(\eta))$. Let $\Delta_\eta = Gal(\mathbf{Q}(\eta)/\mathbf{Q})$. Since $\hat{f}^\sigma = \hat{f}$ for all $\sigma \in \Delta_\eta$, we have an action of $\Delta_\eta$ on the Selmer group $Sel_{O_K, \emptyset}^{\mathbf{Q}(\eta)}(ad^0(\rho_f))$. Moreover the inflation-restriction exact sequence in Galois cohomology gives us the following isomorphism:

$$Sel_{O_K, \emptyset}^{\mathbf{Q}(\eta)}(ad^0(\rho_f))[\eta^{-1}] = Sel_{O_K, \emptyset}(ad^0(\rho_f) \otimes \eta^{-1})$$

Also it follows from Fujiwara's work generalizing Taylor-Wiles' results to totally real fields (cf. [7]) that

$$\Omega^1_{R_{\hat{f}}/O_K} \otimes_{R_{\hat{f}}, \lambda_{\hat{f}}} O_K \cong Sel_{O_K, \emptyset}^{\mathbf{Q}(\eta)}(ad^0(\rho_f))^*$$

where $R_{\hat{f}}$ is the local component of the Hecke algebra of minimal level associated to $\hat{f}$ and $\lambda_{\hat{f}}$ is the character of $R_{\hat{f}}$ given by $\hat{f}$. Since this isomorphism is compatible with the action of $\Delta_\eta$, by $(\star)$, we have that

$$(\star\star) \quad (\Omega^1_{R_{\hat{f}}} \otimes_{R_{\hat{f}}, \lambda_{\hat{f}}} O_K)[\eta^{-1}] \neq 0.$$

This implies that the action of $\Delta$ on the local component $R_{\hat{f}}$ is non trivial and thus $R_{\hat{f}}$ contains irreducible component which are not base change from $\mathbf{Q}$. This proves what we claimed.∎

**Remarks:** a) By a similarly proof, one sees that the same result holds for Hida families of cusp forms.

b) It seems possible to prove the converse of this result when $\chi$ is quadratic. We give a sketch of proof: Assume you have a congruence between $\hat{f}$ and $g$ and consider their theta lifts from $GL_{2/\mathbf{Q}(\eta)}$ to $GSp_{4/\mathbf{Q}}$. Then $\Theta(\hat{f}) \equiv \Theta(g)$ but the eigenvalues of $\Theta(\hat{f})$ are given by the Eisenstein character $\lambda_{Eis(f,\eta)}$ and since $g$ is not a base change lift $\Theta(g)$ should be cuspidal. Therefore $p$ divides the Eisenstein Ideal. We thank D. Prasad for suggesting to us this example of congruences.